\documentclass[preprint,12pt]{elsarticle}

\usepackage[T1]{fontenc}
\usepackage{lmodern}
\usepackage{amsmath,amssymb,amsfonts}
\usepackage{amsthm}
\usepackage{booktabs}
\usepackage{array}
\usepackage{url}
\usepackage[hidelinks]{hyperref}

\newtheorem{theorem}{Theorem}[section]
\newtheorem{proposition}[theorem]{Proposition}
\newtheorem{lemma}[theorem]{Lemma}
\newtheorem{corollary}[theorem]{Corollary}
\newtheorem{maintheorem}{Theorem}

\theoremstyle{definition}
\newtheorem{definition}[theorem]{Definition}
\newtheorem{example}[theorem]{Example}
\theoremstyle{remark}
\newtheorem{remark}[theorem]{Remark}

\newcommand{\M}{\operatorname{M}}
\newcommand{\N}{\mathbb{N}}
\newcommand{\Z}{\mathbb{Z}}
\newcommand{\Units}[1]{(\Z/#1\Z)^{\times}}

\allowdisplaybreaks

\journal{Journal of Number Theory}

\begin{document}

\begin{frontmatter}

\title{The golden ratio and the spacings of multiples modulo \texorpdfstring{$p$}{p}}

\author[1]{Prashant Chaudhary\corref{cor1}}
\ead{prashant@prashantchaudhary.com}
\author[2]{Kapil Khandelwal}
\ead{kapil.khandelwal@nd.edu}
\cortext[cor1]{Corresponding author.}

\affiliation[1]{organization={Independent Researcher},
            city={Santa Clara},
            state={CA},
            country={USA}}
\affiliation[2]{organization={Department of Civil and Environmental Engineering and Earth Sciences, University of Notre Dame},
            city={Notre Dame},
            state={IN},
            country={USA}}

\begin{abstract}
How evenly can the multiples of a unit modulo $p$ spread around $\Z/p\Z$ when every stage counts? For an integer $p\ge2$ and a unit $r$ modulo $p$, place $0,r,2r,\dots$ on $\Z/p\Z$ one at a time, and let $\Phi(p,r)$ be the sum, over all $p$ stages, of $\binom g2$ over the arc lengths $g$. Here $p\,\Phi(p,r)$ counts the empty periodic boxes of the lattice $\{(i,ir)\}$. We show that the minimum $m(p)$ of $\Phi(p,r)$ over the units satisfies $\liminf m(p)/(p^2\log p)=C_{\mathrm g}$, where $C_{\mathrm g}=\tau^2/(10\log\tau)=0.544\ldots$ and $\tau$ is the golden ratio; the Cauchy--Schwarz inequality alone gives the constant $1/2$. The proof is elementary and gives $\Phi(p,r)\ge C_{\mathrm g}p^2\log(p/2)+p^2/8+p/4$ for every $p\ge2$ and every unit $r$, and the ratio $m(p)/(p^2\log p)$ tends to $C_{\mathrm g}$ along the Fibonacci numbers. Up to an error of $O(p^2\log p)$, $\Phi(p,r)$ equals $p^2/6$ times the sum of the partial quotients of $r/p$, so $m(p)=O(p^2\log p)$ if and only if every $p$ has a unit whose partial quotients sum to $O(\log p)$. The liminf is not a limit: an elementary count shows that the ratio exceeds $C_{\mathrm g}$ by a fixed amount for almost all $p$. Computer-assisted certificates show that only $O((\log X)^{69})$ integers $p\le X$ have $m(p)\le C_{\mathrm g}p^2\log p$, and, together with theorems of Kan and of Pollicott and Vytnova, place the ratio between $0.54538$ and $0.584$ for a set of $p$ of density one. These are asymptotic statements: among $2\le p\le3000$ only $p=506$ has $m(p)>C_{\mathrm g}p^2\log p$. For coprime $p$ and $q$, a reciprocity law relates $\Phi(p,q)$ and $\Phi(q,p)$ through the mass of the numerical semigroup $\langle p,q\rangle$, and the three quantities form a binary quadratic form that transforms along the continued fraction of $q/p$.
\end{abstract}


\end{frontmatter}

\section{Introduction}\label{sec:intro}

Step around the circle $\Z/p\Z$ by a fixed unit $r$, marking $0,r,2r,\dots$ one point at a time; here $p\ge2$ is any integer, not necessarily prime. After $h$ marks the circle is cut into $h$ arcs, and by the three-distance theorem the arcs take at most three lengths \cite{Sos1958,Suranyi1958,Swierczkowski1958}. How should $r$ be chosen so that the marks are spread evenly at every stage of the walk? To make the question precise, charge an arc of length $g$ the amount $\binom g2$, so that one long arc costs more than several short ones of the same total. For $1\le h\le p$ let $\phi_{p,r}(h)$ be the total charge at stage $h$, and put
\[
\Phi(p,r)=\sum_{h=1}^{p}\phi_{p,r}(h).
\]
The functional first came up in the companion paper \cite{Undershoot}, whose question arose from filling advertising breaks with ads of a few fixed durations: there, for coprime $0<a<b$, $2\Phi(b,a)-\binom b2$ is the limiting total undershoot of the sumsets of $\{0,a,b,2b\}$ (Remark~\ref{rem:link}). Two facts make $\Phi$ natural. At each stage $\phi_{p,r}(h)=\frac12\sum_ig_i^2-\frac p2$, where $g_1,\dots,g_h$ are the arc lengths, is an affine function of the sum of the squared arcs, a classical statistic of spacings \cite{Greenwood1946,Moran1947}, and $p\,\Phi(p,r)$ is the number of empty periodic boxes of the lattice $\{(i,ir)\}\subset(\Z/p\Z)^2$ (Lemma~\ref{lem:rect}).

For $p=13$ the units $r\equiv\pm1$ give $\Phi=364$, the largest possible value, $r\equiv\pm2$ and $\pm6$ give $239$, $r\equiv\pm3$ and $\pm4$ give $214$, and $r\equiv\pm5$ give the least value, $198$. The best steps $5$ and $8$ are, like $13$, Fibonacci numbers, and for $p=21$ the least value, $629$, is again attained at $r=8$ and $r=13$. The golden ratio $\tau=(1+\sqrt5)/2$ is behind this pattern. We determine the lower limit of the minimum $m(p)=\min_{r\in\Units p}\Phi(p,r)$, show that for almost all $p$ the ratio $m(p)/(p^2\log p)$ stays a fixed distance above it, and determine the maximum and the mean of $\Phi(p,r)$ over the units. Let $F_n$ be the Fibonacci numbers, with $F_0=0$ and $F_1=1$, and put $C_{\mathrm g}=\tau^2/(10\log\tau)\approx0.5440502$.

\begin{maintheorem}[Theorem~\ref{thm:golden}]\label{thm:A}
For $p\ge2$ and every unit $r$ modulo $p$,
\[
\Phi(p,r)\ge C_{\mathrm g}\,p^2\log\frac p2+\frac{p^2}8+\frac p4,
\]
with equality for $p=2$, and $\Phi(F_n,F_{n-1})=C_{\mathrm g}F_n^2\log F_n+O(F_n^2)$ for $n\ge3$. Consequently
\[
\liminf_{p\to\infty}\frac{m(p)}{p^2\log p}=C_{\mathrm g}\qquad\text{and}\qquad\lim_{n\to\infty}\frac{m(F_n)}{F_n^2\log F_n}=C_{\mathrm g}.
\]
\end{maintheorem}

The Cauchy--Schwarz inequality, applied stage by stage, gives $\Phi(p,r)\ge\frac{p^2}2(H_p-1)$ with $H_p=\sum_{h\le p}1/h$ (Theorem~\ref{thm:extremes}), that is, the constant $\frac12$ in place of $C_{\mathrm g}=0.544\ldots$. The sharp coefficient exceeds this stagewise coefficient by $0.0440\ldots$. The proof of Theorem~\ref{thm:A} is elementary and has no computer-assisted step; the code also checks the polynomial identities it uses (Section~\ref{sec:computation}).

\paragraph{Partial quotients} Write $r/p=[0;a_1,\dots,a_K]$ with $a_K\ge2$ and let $\Sigma(p,r)=a_1+\dots+a_K$ be the sum of its partial quotients. Then
\[
\Phi(p,r)=\frac{p^2}6\,\Sigma(p,r)+O(p^2\log p)
\]
with an absolute constant (Theorem~\ref{thm:pq}), so $\Phi$ refines the sum of partial quotients. The error bound has the order of the minimum of $\Phi$, along the Fibonacci numbers and for almost all $p$ (Theorems~\ref{thm:A} and~\ref{thm:typical-upper}), so the two minimum problems can differ: the units that minimize $\Phi(p,r)$ need not minimize $\Sigma(p,r)$, as a computer-assisted example at $p=221688$ shows (Proposition~\ref{prop:counterexample}). The constant $C_{\mathrm g}$ is not visible from the partial quotients: along the Fibonacci numbers $F_n^2\Sigma(F_n,F_{n-1})/6=(1/(6\log\tau)+o(1))\,F_n^2\log F_n$, with $1/(6\log\tau)=0.346\ldots$, and the error term of Theorem~\ref{thm:pq} accounts for the remaining $0.197\ldots$ of $C_{\mathrm g}$. At the coarser scale the two problems agree. With $\Sigma_{\min}(p)=\min_r\Sigma(p,r)$, one has $m(p)=p^2\Sigma_{\min}(p)/6+O(p^2\log p)$, and $\limsup_{p\to\infty}m(p)/(p^2\log p)$ is finite if and only if $\Sigma_{\min}(p)=O(\log p)$, a conjecture attributed to Moser (see \cite[Sect.~1]{Shkredov2026}; Proposition~\ref{prop:order}). The limsup is finite if the recent preprints \cite{Shkredov2026,Zhang2026} on Zaremba's conjecture are correct (Remark~\ref{rem:zhang}).

The mean of $\Phi(p,r)$ over the units is larger than the lower bound of Theorem~\ref{thm:A} by a factor of order $\log p$: it is $p^2\log^2p/\pi^2+O(p^2\log p)$ (Theorem~\ref{thm:average}). A typical unit lies in between, with $\Phi(p,r)$ of order $p^2\log p\log\log p$ (Remark~\ref{rem:random-unit}), and the maximum is $\binom{p+1}3$, attained exactly at $r\equiv\pm1$ (Theorem~\ref{thm:extremes}). Along the Fibonacci numbers $\Phi(F_n,F_{n-1})/(F_n^2\log F_n)$ approaches $C_{\mathrm g}$ from below, with a deficit of about $0.242/\log F_n$ (Remark~\ref{rem:fib-second-order}). Few denominators come close to the liminf.

\begin{maintheorem}[Theorems~\ref{thm:sparse} and~\ref{thm:penalty}]\label{thm:B}
There is an explicit $\varepsilon_0>0$ such that, with $\theta=(1+\log31)/10=0.4433\ldots$,
\[
\#\bigl\{2\le p\le X:\ m(p)\le(C_{\mathrm g}+\varepsilon_0)\,p^2\log p\bigr\}=O(X^{\theta}\log X).
\]
With a computer-assisted certificate,
\begin{align*}
\#\{2\le p\le X:\ m(p)\le C_{\mathrm g}\,p^2\log p\}&=O\bigl((\log X)^{69}\bigr),\\
\#\{2\le p\le X:\ m(p)\le(C_{\mathrm g}+\tfrac1{200000})\,p^2\log p\}&=O(X^{1/90}),\\
\#\{2\le p\le X:\ m(p)\le(C_{\mathrm g}+\tfrac1{750})\,p^2\log p\}&=O(X^{39/40}).
\end{align*}
\end{maintheorem}

One can take $\varepsilon_0=1/10800000$. The first count is elementary, and since $\theta<1$ it gives $m(p)>(C_{\mathrm g}+\varepsilon_0)p^2\log p$ for almost all $p$, so $m(p)/(p^2\log p)$ has no limit. On the upper side, units with bounded partial quotients, the subject of Zaremba's conjecture \cite{Zaremba1972,BourgainKontorovich2014,Huang2015}, give upper bounds for $m(p)$ on a set of density one.

\begin{maintheorem}[Theorem~\ref{thm:typical-upper}; computer-assisted]\label{thm:C}
There is $c_0>0$ such that, for all $p\le X$ outside a set of size $O(Xe^{-c_0\sqrt{\log X}}+X^{1-1/1250})$,
\[
C_{\mathrm g}+\frac1{750}<\frac{m(p)}{p^2\log p}<\frac{73}{125}=0.584 .
\]
In particular these bounds hold for a set of integers of density one and for a set of primes of relative density one.
\end{maintheorem}

Both bounds are computer-assisted: the lower bound is the last count of Theorem~\ref{thm:B}, and the upper bound uses a transfer-matrix certificate together with theorems of Kan \cite{Kan2021} and of Pollicott and Vytnova \cite{PollicottVytnova2022}. For every $p\ge3$ the best upper bound we have is $m(p)\le\frac2{\pi^2}p^2\log p\log\log p+O(p^2\log p)$, from a theorem of Aistleitner, Borda, and Hauke (Remark~\ref{rem:abh-upper}).

The counts are asymptotic and are informative only for very large $X$: the explicit majorant behind the first count of Theorem~\ref{thm:B} exceeds $X$ for $X\le10^{900000}$, and the one behind $O((\log X)^{69})$ for $X\le10^{83}$. At accessible sizes the typical behavior is not yet visible: among $2\le p\le3000$ only $p=506$ has $m(p)>C_{\mathrm g}p^2\log p$ (Remark~\ref{rem:scale}).

\paragraph{Reciprocity} The second argument of $\Phi$ may exceed the first, and the swap $\Phi(p,q)\leftrightarrow\Phi(q,p)$ is governed by a third quantity. For coprime $p,q\ge1$ the numerical semigroup $\langle p,q\rangle=\{up+vq:u,v\in\Z_{\ge0}\}$ has finitely many gaps, the positive integers outside it. Its \emph{mass} $\M(p,q)$, defined in \cite{Undershoot}, is the sum over the gaps $x$ of the distance from $x$ down to the largest element of $\langle p,q\rangle$ below $x$; it is also the genus of a two-dimensional semigroup lifted from $\langle p,q\rangle$ \cite[Cor.~3.7(2), with $d_j=1$]{BhardwajCistoDeshpande2026}.

\begin{maintheorem}[Theorems~\ref{thm:reciprocity} and~\ref{thm:closed-form}]\label{thm:D}
For coprime positive integers $n,a$,
\[
\Phi(n+a,a)=\Phi(n,a)+\Phi(a,n)+2\M(a,n)+\binom{n+a}{2}.
\]
For coprime $q>p\ge1$ the binary quadratic form $Q_{q,p}(x,y)=\Phi(q,p)x^2+2\M(p,q)xy+\Phi(p,q)y^2$ satisfies
\[
Q_{q,p}(x,y)=P(p,s,c)\,x^2+Q_{p,s}(cx+y,\,x)\qquad(q=cp+s,\ 0\le s<p),
\]
where $P(n,s,c)=\sum_{j=1}^c\binom{jn+s}2$ and $Q_{1,0}\equiv0$. Along Euclid's algorithm this writes $\Phi(q,p)$, $\M(p,q)$, and $\Phi(p,q)$ as sums of $O(1+\log p)$ explicit terms.
\end{maintheorem}

The reciprocity law relates the three moduli $n$, $a$, and $n+a$. Integrated along Euclid's algorithm (Lemma~\ref{lem:integrated}) it relates $\Phi(q,p)$ and $\Phi(p,q)$, with the mass as the cross term, in the manner of the Dedekind reciprocity law \cite{RademacherGrosswald1972}, and the step of $Q$ is the matrix $\left(\begin{smallmatrix}c&1\\1&0\end{smallmatrix}\right)$ of the continued fraction. The upper bound of Theorem~\ref{thm:pq} is read off this expansion. Twelve times the Dedekind sum $s(r,p)$ is the alternating sum of the partial quotients of $r/p$ plus a bounded explicit term \cite[Thm.~1]{Hickerson1977}, while $\Phi(p,r)/p^2$ is one sixth of their plain sum only up to $O(\log p)$ (Theorem~\ref{thm:pq}). Nor is $\Phi$ a function of $p$ and $s(r,p)$: $s(8,65)=s(18,65)=0$, while $\Phi(65,8)=12048$ and $\Phi(65,18)=8875$. We know of no exact identity between $\Phi$ and Dedekind sums.

\paragraph{Numerical semigroups} The mass is the semigroup side of $\Phi$. Along a residue class $q\equiv r\pmod p$ the mass of $\langle p,q\rangle$ is affine in $q$ with slope $\Phi(p,r)/p$ (Corollary~\ref{cor:affine}): for $q=18,31,44,57$ the masses of $\langle13,q\rangle$ are $237,435,633,831$, and each step adds $\Phi(13,5)=198$. In the study of concordance invariants of torus knots, Wang \cite[Thm.~1.4 with Lemma~2.4]{Wang2020} found a recurrence for the numbers of gap runs of each length of $\langle p,cp+r\rangle$, through the same blocks of length $p$; weighting it gives the affine law in the form \eqref{eq:mass-only}. The new point here is the identification of the slope with $\Phi$. The companion paper \cite{Undershoot} uses the mass for sumsets; the two papers share the identity of Remark~\ref{rem:link}.

\paragraph{Method} Everything rests on the insertion identity: adding the point $jr$ splits one arc into two, of lengths $d_j^-$ and $d_j^+$, and lowers $\phi_{p,r}$ by $d_j^-d_j^+$. Grouping the stages by the blocks of length $p$ of $\langle p,q\rangle$ gives the Euclidean step (Theorem~\ref{thm:euclid-step}), and the Euclidean step together with a complement identity gives reciprocity.

For Theorem~\ref{thm:A} we run the subtractive Euclidean algorithm on the two arc lengths $r$ and $p-r$ of the second stage. Each subtraction is a step with a state $(x,y)$ of two ratios larger than one. Up to a bounded term, $\Phi(p,r)/p^2$ is the sum of costs $\Gamma(x,y)$ over the steps, while $\log(p/2)$ is exactly the sum of lengths $\Lambda(x,y)$ (Proposition~\ref{prop:action}). The state map fixes $(\tau,\tau)$, where the cost per unit length is $C_{\mathrm g}$. The residual $\Gamma-C_{\mathrm g}\Lambda$ is positive whenever one of the two ratios is at least two; it can be negative only when both ratios are below two, and a three-point inequality between neighboring steps absorbs it there. Summing over the steps gives Theorem~\ref{thm:A}.

For $r\le p/2$ the path has $\sum_i(a_i-1)-1$ no-flip steps (Proposition~\ref{prop:action}), each with residual bounded below; counting the continued fractions with few no-flips gives the first count of Theorem~\ref{thm:B}, and a certified polynomial inequality at the moves between flips and no-flips (\ref{app:penalty}) gives the other three. Theorem~\ref{thm:C} combines an exact formula for $\Phi$ over the digits of a continued fraction with the density-one method of Bourgain and Kontorovich \cite{BourgainKontorovich2014}.

\paragraph{Related work} Linear statistics of the arcs of $\{0,r,\dots,(h-1)r\}$ are classical through the three-distance theorem \cite{vanRavenstein1988,AlessandriBerthe1998}. In the continuous setting, the largest empty box of a lattice $\{(i/p,\{ir/p\})\}$ on the torus is minimized by the Fibonacci lattices \cite{BreneisHinrichs2020}; $p\,\Phi(p,r)$ counts the discrete empty periodic boxes instead. For the sequence $(k\alpha)$ the golden ratio is optimal for several criteria on how new points divide the gaps \cite[Thms.~5--7]{vanRavenstein1989}, and the golden Kronecker sequence has a sharp logarithmic discrepancy constant \cite{Dupain1979}. Among arbitrary sequences of points on the circle the log sequence of de Bruijn and Erd\H{o}s \cite{deBruijnErdos1949} is optimal, in the limit superior over the stages, for a class of monotone measures of the arcs that includes the normalized sum of squares $h\sum_is_i^2$ of the $h$ arc lengths $s_i$ \cite[Thm.~5]{Ramshaw1978}; for $\Phi$, the golden ratio becomes optimal once the points are the multiples of one number (Theorem~\ref{thm:A}). The mean of the sum of partial quotients over all numerators $1\le m\le p$ is due to Yao and Knuth \cite{YaoKnuth1975}; the mean of $\Sigma(p,r)$ over the units was found for prime $p$ by Popov \cite{Popov1981} and for all $p$ by Panov \cite{Panov1980,Panov1982} and by Liehl \cite{Liehl1983}, and we use the form of Panov \cite{Panov1982}. Zaremba's conjecture arose from good lattice points for numerical integration, whose discrepancy is controlled by the partial quotients of $r/p$ \cite{Niederreiter1978}. Its density-one forms are due to Bourgain and Kontorovich \cite{BourgainKontorovich2014} and to Huang \cite{Huang2015}, and Kan \cite{Kan2021} lowers the dimension threshold of their density-one theorem, including the lower bound for the number of numerators that Section~\ref{sec:sparsity} uses. For numerical semigroups, linear gap statistics such as the genus and the Sylvester power sums are read off one Ap\'ery set and are polynomials in $p$ and $q$ for two generators \cite{Selmer1977,BrownShiue1993,Rodseth1994,Tuenter2006}; the mass is quadratic in the run lengths, and for two generators it is, on the integers $q$ coprime to $p$, quasi-polynomial in $q$ with period $p$ but not a polynomial in $p$ and $q$ (Corollary~\ref{cor:affine} and Example~\ref{ex:mass-small}). Apart from \cite{Undershoot}, we have not found the sum $\Phi$ over all stages; the mass occurs in \cite{Wang2020} through its run vector and in \cite{BhardwajCistoDeshpande2026} as a genus, and is studied as a statistic in its own right in \cite{Undershoot}.

\paragraph{Organization} Section~\ref{sec:prelim} fixes notation. Section~\ref{sec:spacings} studies $\phi$ and $\Phi$, Section~\ref{sec:reciprocity} proves Theorem~\ref{thm:D}, and Section~\ref{sec:pq} compares $\Phi$ with the sum of partial quotients and treats its maximum and mean. Sections~\ref{sec:golden} and~\ref{sec:sparsity} prove Theorems~\ref{thm:A} to~\ref{thm:C}, and Section~\ref{sec:minimizers} describes the minimizers. Open problems follow in Section~\ref{sec:open}. Section~\ref{sec:computation} explains how the exact computations are checked and rerun from the code \cite{code2026}, and the appendices hold the long proofs.

\section{Notation}\label{sec:prelim}

\paragraph{Residues and continued fractions} For integers $x$ and $n\ge1$, $\langle x\rangle_n$ is the least nonnegative residue of $x$ modulo $n$, and $\Units p$ is the group of units modulo $p$; a unit $r$ is identified with its representative in $[1,p-1]$. For $1\le r<p$ coprime to $p$ we write the canonical continued fraction
\[
\frac rp=[0;a_1,\dots,a_K],\qquad a_K\ge2,\qquad \Sigma(p,r)=\sum_{i=1}^Ka_i,
\]
and $\Sigma_{\min}(p)=\min_r\Sigma(p,r)$. The digit sum is the same for the expansion ending in $1$, and $\Sigma(p,r)=\Sigma(p,p-r)$. Since the smallest denominator of a canonical expansion of length $K$ is $F_{K+2}\ge\tau^K$, where $F_0=0$, $F_1=1$ are the Fibonacci numbers and $\tau=(1+\sqrt5)/2$,
\begin{equation}\label{eq:length}
K\le\frac{\log p}{\log\tau}.
\end{equation}
The letters $c$ and $c_i$ ($i\ge1$) denote quotients of Euclid's algorithm, while $c_0$ (Theorem~\ref{thm:C}) and $c_{\mathrm F}$ (Remark~\ref{rem:fib-second-order}) are constants. All logarithms are natural, and $C_{\mathrm g}=\tau^2/(10\log\tau)$ as in the introduction.

\paragraph{Numerical semigroups and the mass} We follow \cite{RosalesGarciaSanchez2009}, with $\N=\{0,1,2,\dots\}$. A numerical semigroup $S$ is a submonoid of $(\N,+)$ whose complement in $\N$ is finite; the integers in that complement are its gaps. For a positive $m\in S$, the Ap\'ery set of $S$ with respect to $m$ is the set of least elements of $S$ in the residue classes modulo $m$.

\begin{definition}[Mass {\cite[Def.~2.1]{Undershoot}}]\label{def:mass}
For a numerical semigroup $S$ put
\[
\M(S)=\sum_{x\notin S}\bigl(x-\max\{s\in S: s<x\}\bigr),
\]
the sum over the gaps of the distance to the largest element of $S$ below them. For coprime positive integers $p,q$ we write $\M(p,q)=\M(\langle p,q\rangle)$, so $\M(p,q)=\M(q,p)$ and $\M(1,q)=0$.
\end{definition}

A maximal run of $L$ consecutive gaps contributes $1+2+\dots+L=\binom{L+1}2$ to $\M(S)$.

\section{Spacings of multiples}\label{sec:spacings}

\begin{definition}[{\cite[Sect.~2]{Undershoot}}]\label{def:phi}
Let $n\ge1$ and $a$ be coprime integers. For $1\le h\le n$ let $g_1,\dots,g_h$ be the circular spacings of the $h$ points $0,a,2a,\dots,(h-1)a$ in $\Z/n\Z$, so that $g_i\ge1$ and $g_1+\dots+g_h=n$; for $h=1$ the single point $0$ has one spacing, of length $n$. Put
\[
\phi_{n,a}(h)=\sum_{i=1}^h\binom{g_i}{2},\qquad
\Phi(n,a)=\sum_{h=1}^{n}\phi_{n,a}(h),
\]
and $\phi_{n,a}(h)=0$ for $h>n$.
\end{definition}

Both quantities depend only on $a\bmod n$. Always $\phi_{n,a}(1)=\binom n2$ and $\phi_{n,a}(n)=0$, and $\Phi(1,a)=0$. For $a=1$ the stage $h$ has one spacing $n-h+1$ and $h-1$ unit spacings, so
\begin{equation}\label{eq:phi-one}
\phi_{n,1}(h)=\binom{n-h+1}2,\qquad \Phi(n,1)=\binom{n+1}3 .
\end{equation}

\paragraph{Insertion} For $1\le j\le n-1$ put
\[
d_j^-=\min_{1\le t\le j}\langle ta\rangle_n,\qquad
d_j^+=\min_{1\le t\le j}\langle -ta\rangle_n,\qquad
\lambda_{n,a}(j)=\phi_{n,a}(j)-\phi_{n,a}(j+1).
\]
When several pairs occur we write $d_j^\pm(n,a)$. Inserting the point $ja$ into $\{0,a,\dots,(j-1)a\}$ splits one circular spacing into two parts. Their lengths are $d_j^-$ and $d_j^+$, the distances from $ja$ back to the preceding old point and forward to the following one, because the differences between $ja$ and the old points are the multiples $ta$, $1\le t\le j$. This is the mechanism of the three-distance theorem \cite{Sos1958,Suranyi1958,Swierczkowski1958,AlessandriBerthe1998}: $d_j^-$ and $d_j^+$ are the record distances of the orbit $(ta\bmod n)_t$. Since $\binom{L+R}2-\binom L2-\binom R2=LR$,
\begin{equation}\label{eq:insertion}
\lambda_{n,a}(j)=d_j^-d_j^+,\qquad
\sum_{j=1}^{n-1}\lambda_{n,a}(j)=\binom n2,\qquad
\Phi(n,a)=\sum_{j=1}^{n-1}j\,\lambda_{n,a}(j),
\end{equation}
where the last two identities telescope, using $\phi_{n,a}(1)=\binom n2$ and $\phi_{n,a}(n)=0$ (all sums are empty for $n=1$).

\begin{lemma}[Rectangles]\label{lem:rect}
For $n\ge1$ and $\gcd(n,a)=1$,
\[
\Phi(n,a)=\sum_{W\ge1}\sum_{H\ge1}\Bigl(n-\bigl|a[0,W)+[0,H)\bigr|\Bigr),
\]
where $[0,W)$ is the set of integers $0\le i<W$ and the set $a[0,W)+[0,H)$ is taken in $\Z/n\Z$. Terms with $W\ge n$ or $H\ge n$ vanish. Consequently $n\,\Phi(n,a)$ is the number of empty periodic boxes of the lattice $\{(i,ia)\}\subset(\Z/n\Z)^2$: the quadruples $(x,y,W,H)$ with $x,y\in\Z/n\Z$ and integers $W,H\ge1$ for which no $i\in x+[0,W)$ has $ia\in y+[0,H)$.
\end{lemma}

\begin{proof}
If $W\ge n$ or $H\ge n$ the set $a[0,W)+[0,H)$ is all of $\Z/n\Z$. For $1\le W<n$ let $g_1,\dots,g_W$ be the circular spacings of $X=a[0,W)$. A residue lies in $X+[0,H)$ exactly when its backward distance to the most recent point of $X$ is less than $H$, so the arc from a point of $X$ up to the next point contributes $\min(g_i,H)$ elements. Hence $n-|X+[0,H)|=\sum_i(g_i-H)_+$ and $\sum_{H\ge1}\sum_i(g_i-H)_+=\sum_i\binom{g_i}2=\phi_{n,a}(W)$. Summing over $W$ proves the formula. The box $(x,y,W,H)$ is empty exactly when $y\notin ax+(X-[0,H))$. A residue lies in $X-[0,H)$ exactly when its forward distance to the next point of $X$ is less than $H$, so $|X-[0,H)|=\sum_i\min(g_i,H)=|X+[0,H)|$; for $W\ge n$ both sets are all of $\Z/n\Z$. Hence, for each of the $n$ values of $x$, the number of $y$ with $(x,y,W,H)$ empty is the summand for $(W,H)$.
\end{proof}

\begin{proposition}[Symmetries]\label{prop:sym}
Let $p\ge1$ and $\gcd(p,r)=1$, and let $r^*$ be an inverse of $r$ modulo $p$. Then $\phi_{p,r}(h)=\phi_{p,-r}(h)$ for every $h$, and
\[
\Phi(p,r)=\Phi(p,-r)=\Phi(p,r^*).
\]
\end{proposition}

\begin{proof}
The reflection $x\mapsto-x$ of $\Z/p\Z$ maps the points for $r$ to the points for $-r$ and preserves every spacing. Multiplication by $r^*$ is a bijection of $\Z/p\Z$ mapping $r[0,W)+[0,H)$ to $[0,W)+r^*[0,H)$, so the summand of Lemma~\ref{lem:rect} for $(r,W,H)$ equals the summand for $(r^*,H,W)$. Summing over $W$ and $H$ gives $\Phi(p,r)=\Phi(p,r^*)$.
\end{proof}

Inversion preserves the sum but not the individual values: modulo $7$, $2^*=4$ and
\[
(\phi_{7,2}(h))_{h=1}^{7}=(21,11,5,3,2,1,0),\qquad(\phi_{7,4}(h))_{h=1}^{7}=(21,9,6,4,2,1,0),
\]
both with sum $43$. For $p=13$ and $r=5$,
\[
(\phi_{13,5}(h))_{h=1}^{13}=(78,38,23,17,11,9,7,5,4,3,2,1,0),
\]
with sum $198$; for $r=1$ the terms are $\binom{14-h}2$, with sum $364$.

\section{The Euclidean step and reciprocity}\label{sec:reciprocity}

Throughout this section $\phi_{n,a}$, $\Phi(n,a)$, and $\lambda_{n,a}$ are as in Section~\ref{sec:spacings}, and $\M(p,q)$ is the mass of $\langle p,q\rangle$ (Definition~\ref{def:mass}). In $\phi_{n,a}$, $\Phi(n,a)$, and $\lambda_{n,a}$ the first argument is always the modulus; the second is the step, which may exceed the modulus and matters only modulo it. The mass is symmetric in its two arguments. In Section~\ref{subsec:blocks} the modulus is the generator $p$, and from Section~\ref{subsec:reciprocity} on both $\Phi(q,p)$ and $\Phi(p,q)$ occur. The first lemma writes the mass as a sum of the $\phi_{p,q}(h)$, weighted by how many blocks of length $p$ of the semigroup sit at stage $h$.

\subsection{Blocks and the Euclidean step}\label{subsec:blocks}

\begin{lemma}[Blocks]\label{lem:block}
For all coprime $p,q\ge1$, with $\delta_h(q)=\lfloor hq/p\rfloor-\lfloor(h-1)q/p\rfloor$ (which depends on $p$ as well),
\[
\M(p,q)=\sum_{h=1}^{p-1}\delta_h(q)\,\phi_{p,q}(h)=\sum_{j=1}^{p-1}\Bigl\lfloor\frac{jq}{p}\Bigr\rfloor\lambda_{p,q}(j).
\]
\end{lemma}

\begin{proof}
The case $p=1$ is trivial. For $0\le j<p$ let $x_j\in[0,p)$ satisfy $x_jq\equiv j\pmod p$. The least element of $\langle p,q\rangle$ congruent to $j$ is $x_jq$, since any representation $up+vq$ in this class has $v=x_j+tp$ with $t\ge0$. In the block $\{kp,\dots,(k+1)p-1\}$ the integer $kp+j$ lies in the semigroup if and only if $x_jq\le kp+j$, equivalently $x_jq<(k+1)p$ (the two numbers are congruent modulo $p$), equivalently $x_j<h_k:=\lceil (k+1)p/q\rceil$. Thus if $h_k<p$ the residues present in the block are exactly $0,q,\dots,(h_k-1)q$ modulo $p$, and if $h_k\ge p$ the block has no gaps. Multiples of $p$ lie in the semigroup, so every maximal gap run lies in one block, and the gap runs of the block are the circular spacings $g\ge2$ of the present residues, shortened by one (the run after the last present residue ends at the element $(k+1)p$). A run of length $g-1$ contributes $\binom g2$, so the block contributes $\phi_{p,q}(h_k)$, and
\[
\M(p,q)=\sum_{m\ge1}\phi_{p,q}\bigl(\lceil mp/q\rceil\bigr).
\]
For $1\le h\le p-1$ the integers $m\ge1$ with $\lceil mp/q\rceil=h$ are those with $(h-1)q<mp\le hq$, and there are $\delta_h(q)$ of them. This proves the first formula. The second follows from $\phi_{p,q}(h)=\sum_{j=h}^{p-1}\lambda_{p,q}(j)$ and $\sum_{h=1}^j\delta_h(q)=\lfloor jq/p\rfloor$.
\end{proof}
The same block decomposition underlies the recurrence of Wang \cite[Lemma~5.1]{Wang2020} for the gap runs of $\langle p,q\rangle$.

For $q<p$ some stages have multiplicity zero; for instance $\delta_1(q)=0$.

\begin{theorem}[Euclidean step]\label{thm:euclid-step}
Let $p\ge2$ and $q=cp+r$ with $c\ge0$, $1\le r<p$, and $\gcd(p,r)=1$. Then
\[
\M(p,q)=c\,\Phi(p,r)+\M(r,p).
\]
In particular $\Phi(p,r)=\M(p,p+r)-\M(r,p)$.
\end{theorem}

\begin{proof}
In Lemma~\ref{lem:block}, $\phi_{p,q}=\phi_{p,r}$ and $\delta_h(q)=c+\delta_h(r)$. Since $\phi_{p,r}(p)=0$, the first formula gives $\M(p,q)=c\Phi(p,r)+\M(p,r)$, and $\M(p,r)=\M(r,p)$.
\end{proof}

\begin{corollary}[Slope and intercept]\label{cor:affine}
Let $p\ge2$ and let $1\le r<p$ be coprime to $p$. On the progression $q\equiv r\pmod p$, $q\ge r$, the mass is affine in $q$:
\[
\M(p,q)=\alpha(p,r)\,q+\beta(p,r),
\]
with slope and intercept
\[
\alpha(p,r)=\frac{\Phi(p,r)}{p},\qquad \beta(p,r)=\M(r,p)-\frac{r\,\Phi(p,r)}{p}.
\]
\end{corollary}

\begin{proof}
Substitute $c=(q-r)/p$ in Theorem~\ref{thm:euclid-step}.
\end{proof}

Corollary~\ref{cor:affine} is in the spirit of the parametric Frobenius problem, where the Frobenius number of a semigroup whose generators are linear functions of a parameter is eventually a quasi-polynomial in it \cite[Thm.~1.2]{RouneWoods2015}; here the statistic is the mass and the slope is explicit.

\begin{proposition}[Complements]\label{prop:complement}
Under the hypotheses of Theorem~\ref{thm:euclid-step},
\begin{align}
\M(r,p)+\M(p-r,p)&=\Phi(p,r)-\binom p2,\label{eq:complement}\\
\M(p,q)&=(c+1)\M(r,p)+c\,\M(p-r,p)+c\binom p2.\label{eq:mass-only}
\end{align}
\end{proposition}

\begin{proof}
For \eqref{eq:complement} apply Lemma~\ref{lem:block} with modulus $p$ to the smaller generators $r$ and $p-r$; by Proposition~\ref{prop:sym} both use $\phi_{p,r}$. For $1\le h<p$ the number $hr/p$ is not an integer, so $\lfloor hr/p\rfloor+\lfloor h(p-r)/p\rfloor=h-1$; taking differences, $\delta_1(r)+\delta_1(p-r)=0$ and $\delta_h(r)+\delta_h(p-r)=1$ for $2\le h<p$. Adding the two block formulas gives $\sum_{h=2}^{p-1}\phi_{p,r}(h)=\Phi(p,r)-\binom p2$. Substituting \eqref{eq:complement} into Theorem~\ref{thm:euclid-step} gives \eqref{eq:mass-only}. Identity \eqref{eq:mass-only} is also what the recurrence of Wang \cite[Thm.~1.4 with Lemma~2.4]{Wang2020} for the numbers of gap runs of each length gives, after weighting a run of $j$ gaps by $\binom{j+1}2$.
\end{proof}

Formula \eqref{eq:mass-only} expresses the mass through the masses $\M(r,p)$ and $\M(p-r,p)$, whose larger generator is the smaller modulus $p$, without reference to $\Phi$. Iterated, it gives a branching recursion which, unlike Theorem~\ref{thm:closed-form} below, may visit on the order of $p$ pairs (for example when $r=1$).

\begin{remark}\label{rem:link}
For coprime $0<a<b$, Theorem~\ref{thm:euclid-step} with $2b-a=b+(b-a)$, Proposition~\ref{prop:sym}, and \eqref{eq:complement} give
\[
\M(a,b)+\M(b,2b-a)=2\Phi(b,a)-\binom b2 .
\]
This identity is \cite[Cor.~5.5]{Undershoot}, where the left side is the limiting total undershoot of the sumsets of $\{0,a,b,2b\}$ and is computed as the sum of a spectrum \cite[Lemma~3.3]{Undershoot}. The question of \cite{Undershoot} comes from filling advertising breaks with ads of a few fixed durations, and that setting is studied in \cite{Pods}.
\end{remark}

\begin{example}\label{ex:mass-small}
(i) For $p=2$, $\M(2,q)=(q-1)/2$: the gaps $1,3,\dots,q-2$ are isolated. (ii) By \eqref{eq:phi-one} and $\M(p,1)=0$, $\M(p,q)=\frac{q-1}{p}\binom{p+1}{3}$ for $q\equiv1\pmod p$. (iii) For $p=3$, $\Phi(3,1)=\Phi(3,2)=4$ and
\[
\M(3,q)=\begin{cases}(4q-4)/3,& q\equiv1\pmod3,\\ (4q-5)/3,& q\equiv2\pmod3.\end{cases}
\]
Hence the mass is not a polynomial in the generators: a polynomial in $(p,q)$ restricted to $p=3$ is a polynomial in $q$, and it cannot agree with two different affine functions on two infinite progressions.
\end{example}

\subsection{Reciprocity}\label{subsec:reciprocity}

The second index of $\Phi$ may exceed the first; for instance $\Phi(n+a,a)$ concerns the multiples of $a$ modulo $n+a$. Reciprocity relates it to the two smaller moduli.

\begin{theorem}[Reciprocity]\label{thm:reciprocity}
For coprime positive integers $n,a$,
\[
\Phi(n+a,a)=\Phi(n,a)+\Phi(a,n)+2\M(a,n)+\binom{n+a}{2}.
\]
\end{theorem}

\begin{proof}
Put $N=n+a$, so $N\ge2$, $1\le a<N$, $\gcd(N,a)=1$, and $N-a=n$. By \eqref{eq:complement} with modulus $N$ and $r=a$,
\[
\Phi(N,a)=\binom N2+\M(a,N)+\M(n,N).
\]
We show $\M(a,N)=\M(a,n)+\Phi(a,n)$. If $a=1$ both sides vanish. If $a\ge2$, write $n=ca+s$ with $c\ge0$ and $1\le s<a$; then $N=(c+1)a+s$, and Theorem~\ref{thm:euclid-step} with modulus $a$ gives $\M(a,N)=(c+1)\Phi(a,s)+\M(s,a)$ and $\M(a,n)=c\,\Phi(a,s)+\M(s,a)$, while $\Phi(a,s)=\Phi(a,n)$. In the same way $\M(n,N)=\M(n,a)+\Phi(n,a)$. Substituting and using $\M(n,a)=\M(a,n)$ gives the theorem.
\end{proof}

\subsection{The quadratic form along the continued fraction}

For integers $n\ge1$, $s\ge0$, $c\ge0$ put
\[
\begin{aligned}
P(n,s,c)&=\sum_{j=1}^{c}\binom{jn+s}{2}\\
&=\frac{2n^2c^3+3n(n+2s-1)c^2+(n^2+6ns-3n+6s^2-6s)c}{12},
\end{aligned}
\]
the closed form following from $\binom{jn+s}2=\bigl(n^2j^2+n(2s-1)j+s(s-1)\bigr)/2$.

\begin{lemma}[Integrated reciprocity]\label{lem:integrated}
Let $q>p\ge1$ be coprime and $q=cp+s$ with $0\le s<p$. If $s>0$ then
\[
\Phi(q,p)=\Phi(s,p)+P(p,s,c)+c^2\Phi(p,s)+2c\,\M(s,p).
\]
If $s=0$ then $p=1$ and $\Phi(q,1)=\binom{q+1}3=P(1,0,q)$.
\end{lemma}

\begin{proof}
Let $s>0$, so $p\ge2$. For $0\le j<c$ apply Theorem~\ref{thm:reciprocity} to the coprime pair $(n,a)=(jp+s,p)$. By periodicity $\Phi(p,jp+s)=\Phi(p,s)$, and Theorem~\ref{thm:euclid-step} gives $\M(p,jp+s)=j\Phi(p,s)+\M(s,p)$ (for $j=0$ this is the symmetry of $\M$). Hence
\[
\Phi((j+1)p+s,p)-\Phi(jp+s,p)=(2j+1)\Phi(p,s)+2\M(s,p)+\binom{(j+1)p+s}2,
\]
and summing over $0\le j<c$, with $\sum_{j<c}(2j+1)=c^2$, proves the claim. If $s=0$, coprimality forces $p=1$; for $1\le h\le q$ the points $0,1,\dots,h-1$ have $h-1$ unit spacings and one spacing $q-h+1$, so $\Phi(q,1)=\sum_{h=1}^q\binom{q-h+1}2=\binom{q+1}3$. This case is computed directly; the derivation above does not apply, since Theorem~\ref{thm:euclid-step} requires $p\ge2$.
\end{proof}

\begin{theorem}[Quadratic form]\label{thm:closed-form}
For coprime $q>p\ge1$ let $Q_{q,p}(x,y)=\Phi(q,p)x^2+2\M(p,q)xy+\Phi(p,q)y^2$, and let $Q_{1,0}\equiv0$.
\begin{enumerate}\renewcommand{\labelenumi}{(\roman{enumi})}
\item If $q=cp+s$ with $0\le s<p$, then
\begin{equation}\label{eq:form-step-main}
Q_{q,p}(x,y)=P(p,s,c)\,x^2+Q_{p,s}(cx+y,\,x).
\end{equation}
\item Run Euclid's algorithm through the final zero remainder,
\[
n_0=q,\quad n_1=p,\quad n_{i-1}=c_in_i+n_{i+1}\ (1\le i\le K),\quad n_{K+1}=0,
\]
so that $n_K=1$ and the last quotient $c_K=n_{K-1}\ge2$ is included. Put $u_0=0$, $u_1=1$, $v_0=1$, $v_1=0$, $w_{i+1}=c_iw_i+w_{i-1}$ for $w\in\{u,v\}$, and $P_i=P(n_i,n_{i+1},c_i)$. Then $Q_{q,p}(x,y)=\sum_{i=1}^KP_i(u_ix+v_iy)^2$; that is,
\[
\Phi(q,p)=\sum_{i=1}^KP_iu_i^2,\qquad \M(p,q)=\sum_{i=1}^KP_iu_iv_i,\qquad \Phi(p,q)=\sum_{i=1}^KP_iv_i^2 .
\]
\end{enumerate}
\end{theorem}

The step \eqref{eq:form-step-main} is the matrix $\left(\begin{smallmatrix}c&1\\1&0\end{smallmatrix}\right)$ of the continued fraction step acting on the form, and part (ii) iterates it; the proofs are in \ref{app:closed-form}, and the last step uses the terminal case of Lemma~\ref{lem:integrated}.

In the notation of Section~\ref{sec:prelim}, $K$ is the length of the canonical expansion $p/q=[0;c_1,\dots,c_K]$. The number of summands satisfies $K\le1+\log p/\log\tau$, whatever $q$ is: if $p=1$ then $K=1$, and otherwise $c_2,\dots,c_K$ are the digits of the canonical expansion of $s/p$, to which \eqref{eq:length} applies. Each $P_i$ is a fixed polynomial expression, so Theorem~\ref{thm:closed-form} evaluates $\Phi$ and $\M$ with $O(1+\log p)$ arithmetic operations (counted without regard to bit size). For example $16/7=[2;3,2]$ has $(P_1,P_2,P_3)=(156,34,1)$, $(u_i)=(1,2,7)$, $(v_i)=(0,1,3)$, whence $\Phi(16,7)=341$, $\M(7,16)=89$, and $\Phi(7,16)=43$.

\begin{corollary}\label{cor:slope-intercept}
For coprime $1\le p<q$, with $(n_i)$ and $(c_i)$ the chain of Theorem~\ref{thm:closed-form},
\[
\M(p,q)=\sum_{i=1}^{K-1}c_i\,\Phi(n_i,n_{i+1}).
\]
\end{corollary}

\begin{proof}
Apply Theorem~\ref{thm:euclid-step} along the chain; it ends at $\M(1,n_{K-1})=0$.
\end{proof}

\begin{corollary}\label{cor:cauchy}
For coprime positive integers $p$ and $q$,
\[
\M(p,q)^2\le\Phi(q,p)\,\Phi(p,q).
\]
\end{corollary}

\begin{proof}
Both sides are symmetric in $p$ and $q$, and for $p=q=1$ they vanish, so let $q>p$. Every $P_i$ is a sum of binomial coefficients, hence nonnegative, so by Theorem~\ref{thm:closed-form} the form $Q_{q,p}$ is a sum of nonnegative multiples of squares of linear forms. It is therefore positive semidefinite, and its discriminant $4\M(p,q)^2-4\Phi(q,p)\Phi(p,q)$ is at most zero.
\end{proof}

\section{Partial quotients, the maximum, and the mean}\label{sec:pq}

This section compares $\Phi(p,r)$ with the sum of partial quotients $\Sigma(p,r)$ and describes the maximum and the mean of $\Phi(p,r)$ over the units $r$ modulo $p$.

\subsection{The maximum and the sum of partial quotients}

\begin{theorem}\label{thm:extremes}
For $p\ge2$ and every unit $r$ modulo $p$, with $H_p=\sum_{h=1}^p1/h$,
\[
\frac{p^2}2\,(H_p-1)\le\Phi(p,r)\le\binom{p+1}3,
\]
and the upper bound is attained if and only if $r\equiv\pm1\pmod p$. Consequently the slope $\alpha(p,r)$ of Corollary~\ref{cor:affine} satisfies $\alpha(p,r)\le(p^2-1)/6$, and $\Phi(p,r)\ge\frac{p^2}2(\log p+\gamma_{\mathrm E}-1)$, where $\gamma_{\mathrm E}$ is Euler's constant.
\end{theorem}

\begin{proof}
At stage $h$ the spacings $g_1,\dots,g_h$ are positive with sum $p$. With $t_i=g_i-1\ge0$, $\sum t_i=p-h$, we have $\phi_{p,r}(h)=\frac12(\sum t_i^2+\sum t_i)\le\frac12((p-h)^2+p-h)=\binom{p-h+1}2$, with equality if and only if at most one $t_i$ is nonzero. Summing over $h$ gives the upper bound, attained for $r=\pm1$ by \eqref{eq:phi-one} and Proposition~\ref{prop:sym}. Conversely, equality forces equality at $h=2$, where the spacings are $r$ and $p-r$, so one of them is $1$. By Cauchy--Schwarz, $\phi_{p,r}(h)=\frac12(\sum g_i^2-p)\ge\frac12(p^2/h-p)$, and summing gives the lower bound. Finally $H_p-\log p$ decreases to $\gamma_{\mathrm E}$, since $\log(1+1/p)>1/(p+1)$.
\end{proof}

\begin{theorem}\label{thm:pq}
For $p\ge2$ and every unit $r$ modulo $p$,
\[
-p^2\min\Bigl\{\frac K6,\ \frac K{12}+\frac{59}{50}\Bigr\}\le\Phi(p,r)-\frac{p^2\Sigma(p,r)}6\le\frac{121}{288}Kp^2 .
\]
In particular $-Kp^2/6\le\Phi(p,r)-p^2\Sigma(p,r)/6\le\frac{121}{288}Kp^2$, and $|\Phi(p,r)-p^2\Sigma(p,r)/6|\le\frac{121}{288\log\tau}\,p^2\log p$.
\end{theorem}

The proof, given in \ref{app:pq}, bounds each summand of Theorem~\ref{thm:closed-form} for the upper side and the refined lower side, and uses a second decomposition of $\Phi$ along the record distances of the orbit $(jr\bmod p)_j$ for the lower bound $-Kp^2/6$. The bound $-Kp^2/6$ is the better one exactly when $K\le14$. In words, $\Phi(p,r)$ is $p^2/6$ times the sum of the partial quotients of $r/p$, up to $O(p^2\log p)$.

\begin{remark}\label{rem:pq-range}
The constants of Theorem~\ref{thm:pq} are not sharp. Over all units $r$ modulo $p$ for $2\le p\le600$, the ratio $(\Phi(p,r)-p^2\Sigma(p,r)/6)/(Kp^2)$ lies in $[-\frac1{12},\,0.1409]$: the lower end is attained only at $p=2$, and the largest value, $0.14087\ldots$, at $(p,r)=(592,81)$ and $(592,95)$. The proved constants are $-\frac16$ and $\frac{121}{288}\approx0.420$.
\end{remark}

\subsection{The mean over units}

Write $\overline\Phi(p)=|\Units p|^{-1}\sum_{r\in\Units p}\Phi(p,r)$ for the mean of $\Phi$, so that the mean slope of Corollary~\ref{cor:affine} is $\overline\Phi(p)/p$.

\begin{theorem}\label{thm:average}
As $p\to\infty$ through all integers,
\[
\overline\Phi(p)=\frac{p^2}{\pi^2}(\log p)^2+O(p^2\log p).
\]
\end{theorem}

\begin{proof}
Let $H(p)=\sum_{r\in\Units p}\Sigma(p,r)$. By Theorem~\ref{thm:pq} and \eqref{eq:length}, $\overline\Phi(p)=\frac{p^2}6H(p)/|\Units p|+O(p^2\log p)$ with an absolute constant. Panov \cite[Thm.~2, eq.~(10)]{Panov1982} (announced in \cite{Panov1980}) proves, for the reduced fractions $k/N$, $(k,N)=1$, with fixed denominator $N$ (the digit sum is the same for either continued fraction expansion),
\[
\begin{split}H(N)&=\frac6{\pi^2}|\Units N|\log^2N+A_1|\Units N|\log N+A_2|\Units N|Q(N)\\&\quad+A_3|\Units N|+O_\varepsilon(N^{7/8+\varepsilon}),\end{split}
\]
with absolute constants $A_1,A_2,A_3$ and $Q(N)=O((\log\log N)^2)$. For fixed $\eta>0$, $p/|\Units p|=\prod_{\ell\mid p}\ell/(\ell-1)\le C_\eta p^\eta$, because $\ell/(\ell-1)\le\ell^\eta$ for all large primes $\ell$ and $\prod_{\ell\mid p}\ell\le p$. With $\varepsilon=1/16$ and $\eta=1/32$ the remainder divided by $|\Units p|$ is $O(p^{-1/32})$, so $H(p)/|\Units p|=\frac6{\pi^2}\log^2p+O(\log p)$, and the theorem follows.
\end{proof}

\begin{remark}\label{rem:random-unit}
The mean is not the typical value. Aistleitner, Borda, and Hauke \cite[Thm.~1]{AistleitnerBordaHauke2024} show that for $p\ge3$ and $r$ uniform in $\Units p$ the proportion of $r$ with $|\Sigma(p,r)-\frac{12}{\pi^2}\log p\log\log p|\ge t\log p$ is $O_C(1/t)$ for $0<t\le(\log p)^C$. Taking $t=\varepsilon\log\log p$ and using Theorem~\ref{thm:pq}, $\Phi(p,R_p)/(p^2\log p\log\log p)\to2/\pi^2$ in probability when $R_p$ is uniform in $\Units p$.
\end{remark}

\section{The golden-ratio liminf}\label{sec:golden}

The minimum is controlled by an exact formula for $\Phi$ along the subtractive Euclidean algorithm. Let $1\le r<p$ be coprime. Start with lengths and counts $(A,B,U,V)=(r,p-r,1,1)$, interchanging $(A,U)$ with $(B,V)$ if $A<B$. As long as $(A,B)\ne(1,1)$ the current quadruple is a \emph{step} of the path (a step is not one of the stages $h$ of Definition~\ref{def:phi}); its \emph{state} is
\[
(x,y)=\Bigl(\frac AB,\ \frac{U+V}V\Bigr),\qquad x,y>1,
\]
and the algorithm replaces $(A,U)$ by $(A-B,U+V)$ and then interchanges the pairs if $A<B$. The invariant $AV+BU=p$ holds throughout, the lengths stay coprime, and the algorithm stops at $(A,B)=(1,1)$ with counts $(U_*,V_*)$, $U_*+V_*=p$. A step is a \emph{flip} if $1<x<2$ (an interchange follows) and a \emph{no-flip} if $x\ge2$. The states evolve by
\[
\psi_{\mathrm f}(x,y)=\Bigl(\frac1{x-1},\,1+\frac1y\Bigr)\ (1<x<2),\qquad \psi_{\mathrm n}(x,y)=(x-1,\,y+1)\ (x\ge2),
\]
and a no-flip step with $x=2$ is the last step: its subtraction produces $(A,B)=(1,1)$. The path depends only on $\{r,p-r\}$; for $p=2$ it is empty. Put
\[
\Gamma(x,y)=\frac{x^2+y^2}{4(x+y-1)^2},\qquad
\Lambda(x,y)=\frac12\Bigl(\log\Bigl(1+\frac1x\Bigr)+\log\Bigl(1+\frac1y\Bigr)\Bigr).
\]
By Proposition~\ref{prop:action} below, $\Gamma$ is the cost of a step and $\Lambda$ its length: $\Phi(p,r)/p^2$ is the total cost of the path up to a bounded boundary term, and $\log(p/2)$ is its total length. The flip map $\psi_{\mathrm f}$ has the fixed point $(\tau,\tau)$, where $\Gamma(\tau,\tau)=\tau^2/10$ and $\Lambda(\tau,\tau)=\log(1+1/\tau)=\log\tau$; the ratio of cost to length there is $C_{\mathrm g}=\tau^2/(10\log\tau)$, which is why the golden ratio appears. For example, for $(p,r)=(8,5)$ the path is
\begin{center}
\begin{tabular}{@{}ccccc@{}}
\toprule
$(A,B,U,V)$ & $(x,y)$ & type & $\Gamma$ & $\Lambda$ \\
\midrule
$(5,3,1,1)$ & $(5/3,\,2)$ & flip & $61/256$ & $0.4377$ \\
$(3,2,1,2)$ & $(3/2,\,3/2)$ & flip & $9/32$ & $0.5108$ \\
$(2,1,2,3)$ & $(2,\,5/3)$ & no-flip & $61/256$ & $0.4377$ \\
\bottomrule
\end{tabular}
\end{center}
with terminal counts $(U_*,V_*)=(5,3)$: the costs add up to $97/128$, the boundary term of Proposition~\ref{prop:action} is $21/128$, and $\Phi(8,5)=64\cdot\frac{118}{128}=59$, while the lengths add up to $\log4$.

\begin{proposition}\label{prop:action}
Let $(x_j,y_j)$ be the states of the steps of a unit $r$ modulo $p\ge2$. Then
\[
\frac{\Phi(p,r)}{p^2}=\sum_j\Gamma(x_j,y_j)+\frac{p^2+p-U_*V_*-r(p-r)}{4p^2},\qquad \sum_j\Lambda(x_j,y_j)=\log\frac p2 ,
\]
and the boundary term lies in $\bigl[\frac18+\frac1{4p},\ \frac14-\frac1{4p}+\frac1{2p^2}\bigr]$. For every unit $r$ the path has $\Sigma(p,r)-2$ steps (the path and $\Sigma(p,r)$ are unchanged under $r\mapsto p-r$). If moreover $r\le p/2$ and $p/r=[a_1;a_2,\dots,a_K]$ is canonical, then $E-1$ of the steps are no-flips, where $E=\sum_i(a_i-1)$.
\end{proposition}

The proof is in \ref{app:action}. The first identity averages an exact computation along the path for $r$ with the same computation for the reversed path, which is the path of $r^{-1}$; Proposition~\ref{prop:sym} makes the two totals equal. The second telescopes, since a subtraction changes the length sum from $A+B$ to $A$ and the count sum from $U+V$ to $U+2V$.

Let $\varrho=\Gamma-C_{\mathrm g}\Lambda$ be the residual. It takes negative values, but only at states with $1<x,y<2$; for example $\Gamma(\frac{17}{10},\frac{17}{10})=\frac{289}{1152}$ and $\Lambda(\frac{17}{10},\frac{17}{10})=\log\frac{27}{17}$, with ratio $0.5422\ldots<C_{\mathrm g}$. Two local inequalities, proved in \ref{app:golden}, control it: a pointwise bound when $\max(x,y)\ge2$, and a bound for a sum over three consecutive states when $1<x,y<2$.

\begin{lemma}\label{lem:outside}
If $x,y\ge1$ and $\max(x,y)\ge2$, then $\Gamma(x,y)>\kappa\Lambda(x,y)$ with $\kappa=5441/10000>C_{\mathrm g}$. In particular $\varrho(x,y)>0$.
\end{lemma}

\begin{lemma}\label{lem:stencil}
For $1<x,y<2$,
\[
\varrho\bigl(\psi_{\mathrm f}^{-1}(x,y)\bigr)+3\varrho(x,y)+\varrho\bigl(\psi_{\mathrm f}(x,y)\bigr)\ge0,\qquad \psi_{\mathrm f}^{-1}(x,y)=\Bigl(1+\frac1x,\frac1{y-1}\Bigr),
\]
with equality only at $(x,y)=(\tau,\tau)$.
\end{lemma}

\begin{lemma}\label{lem:discharge}
For every unit $r$ modulo $p\ge2$, the residuals $\varrho_j=\varrho(x_j,y_j)$ of its steps satisfy $\sum_j\varrho_j\ge\frac13\sum_{\varrho_j\ge0}\varrho_j\ge0$.
\end{lemma}

\begin{proof}
For $p=2$ the path is empty and both sums vanish; let $p>2$. The first step has $y=(1+1)/1=2$ and the last has $x=2$, so by Lemma~\ref{lem:outside} neither has negative residual. If $\varrho_j<0$, Lemma~\ref{lem:outside} forces $1<x_j,y_j<2$; then step $j$ is a flip, and so is step $j-1$, since a no-flip produces $y>2$. Hence both neighbors exist and are $\psi_{\mathrm f}^{\mp1}(x_j,y_j)$, and Lemma~\ref{lem:stencil} gives $\varrho_{j-1}+3\varrho_j+\varrho_{j+1}\ge0$. Sum over the indices with $\varrho_j<0$, and let $N\le0$ and $P\ge0$ be the sums of the negative and of the nonnegative residuals. Negative neighbors contribute nonpositively, and a nonnegative residual is a neighbor of at most two indices, so $0\le3N+2P$. Hence $P+N\ge P/3$.
\end{proof}

\begin{theorem}\label{thm:golden}
For $p\ge2$ and every unit $r$ modulo $p$,
\[
\Phi(p,r)\ge C_{\mathrm g}\,p^2\log\frac p2+\frac{p^2}8+\frac p4,
\]
with equality for $p=2$. Moreover, for $n\ge3$, $\Phi(F_n,F_{n-1})=C_{\mathrm g}F_n^2\log F_n+O(F_n^2)$. Consequently
\[
\liminf_{p\to\infty}\frac{m(p)}{p^2\log p}=C_{\mathrm g},\qquad\text{and}\qquad \lim_{n\to\infty}\frac{m(F_n)}{F_n^2\log F_n}=C_{\mathrm g}.
\]
\end{theorem}

\begin{proof}
By Proposition~\ref{prop:action}, $\Phi(p,r)/p^2=C_{\mathrm g}\log(p/2)+\sum_j\varrho_j+B$ with $B\ge\frac18+\frac1{4p}$, and $\sum_j\varrho_j\ge0$ by Lemma~\ref{lem:discharge}. For $p=2$ both sides equal $1$.

Let $n\ge4$, $p=F_n$, $r=F_{n-1}$. The path has the $n-3$ states $x_j=F_{n-j-1}/F_{n-j-2}$, $y_j=F_{j+3}/F_{j+2}$, $0\le j\le n-4$: the lengths run down the Fibonacci numbers, every step but the last is a flip, and flips act on counts by $(U,V)\mapsto(V,U+V)$. Each ratio $F_{k+1}/F_k$ with $2\le k\le n-2$ occurs once among the $x_j$ and once among the $y_j$, and these ratios lie in $[3/2,2]$. Binet's formula gives $F_{k+1}-\tau F_k=(-1/\tau)^k$, and $F_k\ge\tau^{k-2}$, so $|F_{k+1}/F_k-\tau|\le\tau^{2-2k}$. On $[3/2,2]^2$,
\[
\partial_x\Gamma(x,y)=\frac{xy-x-y^2}{2(x+y-1)^3},
\]
whose numerator lies in $[-\frac52,-\frac54]$ and whose denominator is at least $16$; so $|\partial_x\Gamma|\le\frac5{32}$, and $|\partial_y\Gamma|\le\frac5{32}$ by symmetry. With $\Gamma(\tau,\tau)=2\tau^2/(4(2\tau-1)^2)=\tau^2/10$, because $(2\tau-1)^2=5$,
\begin{equation}\label{eq:fib-lipschitz}
\Bigl|\sum_j\Gamma(x_j,y_j)-(n-3)\frac{\tau^2}{10}\Bigr|\le\frac5{32}\cdot2\sum_{k\ge2}\tau^{2-2k}=\frac5{16\tau}<\frac15 .
\end{equation}
The boundary term is bounded and $C_{\mathrm g}\log F_n=C_{\mathrm g}n\log\tau+O(1)=(n-3)\tau^2/10+O(1)$, so $\Phi(F_n,F_{n-1})=C_{\mathrm g}F_n^2\log F_n+O(F_n^2)$. Together with the lower bound, which gives $m(p)/(p^2\log p)\ge C_{\mathrm g}+O(1/\log p)$ for all $p$, this proves both limits.
\end{proof}

\begin{remark}[The approach from below]\label{rem:fib-second-order}
The limit
\[
c_{\mathrm F}=\lim_{n\to\infty}\Bigl(\frac{\Phi(F_n,F_{n-1})}{F_n^2}-C_{\mathrm g}\log F_n\Bigr)
\]
exists and is negative, so the Fibonacci ratios approach $C_{\mathrm g}$ from below, at rate $1/\log p$:
\[
\frac{\Phi(F_n,F_{n-1})}{F_n^2\log F_n}=C_{\mathrm g}+\frac{c_{\mathrm F}+o(1)}{\log F_n},\qquad c_{\mathrm F}\approx-0.24232 .
\]
Indeed, here $U_*V_*=r(p-r)=F_{n-1}F_{n-2}$, so the boundary term of Proposition~\ref{prop:action} tends to $(1-2\tau^{-3})/4$, and $\log F_n-n\log\tau\to-\log\sqrt5$. Writing $\Gamma(x,y)-\Gamma(\tau,\tau)$ as $[\Gamma(x,\tau)-\Gamma(\tau,\tau)]+[\Gamma(\tau,y)-\Gamma(\tau,\tau)]$ plus a term bounded by a constant times $|x-\tau||y-\tau|$, the sum in \eqref{eq:fib-lipschitz} converges: the first two parts are partial sums of absolutely convergent series in $k$, and the cross terms add up to $O(n\tau^{-2n})$. Since $\Gamma$ is symmetric, this gives
\[
c_{\mathrm F}=2\sum_{k\ge2}\Bigl(\Gamma\Bigl(\frac{F_{k+1}}{F_k},\tau\Bigr)-\frac{\tau^2}{10}\Bigr)-\frac{3\tau^2}{10}+C_{\mathrm g}\log\sqrt5+\frac{1-2\tau^{-3}}4,
\]
a series whose terms decay like $\tau^{-2k}$. Evaluating it with rational bounds gives the digits above (Section~\ref{sec:computation}); this evaluation is not used in any proof. By \eqref{eq:fib-lipschitz},
\[
c_{\mathrm F}\le-\frac{3\tau^2}{10}+C_{\mathrm g}\log\sqrt5+\frac{1-2\tau^{-3}}4+\frac5{16\tau}<-0.022,
\]
while the lower bound of Theorem~\ref{thm:golden} gives $c_{\mathrm F}\ge\frac18-C_{\mathrm g}\log2=-0.2521\ldots$; the second-order term of that bound is thus within $0.01$ of the Fibonacci value. The ratio is $0.526\ldots$ at $F_{30}=832040$ and $0.538\ldots$ at $F_{90}$, still about $0.006$ below $C_{\mathrm g}$.
\end{remark}

Golden-ratio and Fibonacci extremality is classical for gap-division criteria of the sequence $(k\alpha)$ \cite[Thms.~5--7]{vanRavenstein1989} and for the largest empty box of lattices on the torus \cite{BreneisHinrichs2020}, and sharp logarithmic constants are known in related problems: the star discrepancy $D_N^*$ of the golden Kronecker sequence satisfies $\limsup_N ND_N^*/\log N=3/(20\log\tau)$ \cite{Dupain1979}, and the sum of partial quotients satisfies $p\le F_{\Sigma(p,r)+1}$ for every $p$ and every unit $r$ (splitting a digit $a\ge2$ into $1,a-1$ does not decrease the continuant), with equality along the Fibonacci numbers, so that $\liminf\Sigma_{\min}(p)/\log p=1/\log\tau$. For the periodic $L_2$ discrepancy of plane point sets, whose square is at least of order $\log N$ for every set of $N$ points, the Fibonacci lattices give an explicit upper bound for the optimal constant \cite[Sect.~3]{HinrichsKritzingerPillichshammer2021}, and for the symmetrized Fibonacci sets an exact formula for the $L_2$ discrepancy gives a numerically evaluated constant \cite[Thm.~2.11]{BilykTemlyakovYu2012}; we do not know of a matching lower bound in either case. We have not found an earlier determination of the sharp liminf for a stage-summed quadratic spacing functional of the multiples of a residue. Whether $m(F_n)=\Phi(F_n,F_{n-1})$ for all $n\ge3$ is open (Section~\ref{sec:open}); it holds for $3\le n\le28$ (Section~\ref{sec:computation}).

\section{Sparsity and the typical size of the minimum}\label{sec:sparsity}

\subsection{Sparsity of near-golden denominators}

Lemma~\ref{lem:discharge} keeps one third of the positive residuals, and every no-flip step has a residual bounded below by an absolute constant. This gives an elementary penalty for digits larger than one.

\begin{lemma}\label{lem:noflip}
At a no-flip state, $\varrho(x,y)>1/360000$. Consequently, if $N_0$ is the number of no-flip steps of a unit $r$ modulo $p\ge2$, then
\[
\frac{\Phi(p,r)}{p^2}\ge C_{\mathrm g}\log\frac p2+\frac18+\frac1{4p}+\frac{N_0}{1080000}.
\]
\end{lemma}

\begin{theorem}\label{thm:sparse}
Let $\varepsilon_0=1/10800000$. Then
\[
\#\bigl\{2\le p\le X:\ m(p)\le(C_{\mathrm g}+\varepsilon_0)\,p^2\log p\bigr\}=O(X^{\theta}\log X),
\]
where $\theta=(1+\log31)/10=0.4433\ldots<\frac12$. Hence $m(p)>(C_{\mathrm g}+\varepsilon_0)p^2\log p$ for a set of $p$ of density one,
\[
\limsup_{p\to\infty}\frac{m(p)}{p^2\log p}\ge C_{\mathrm g}+\varepsilon_0>\liminf_{p\to\infty}\frac{m(p)}{p^2\log p},
\]
and there is no constant $A$ with $m(p)\sim Ap^2\log p$ along all integers $p$.
\end{theorem}

The proofs of Lemma~\ref{lem:noflip} and Theorem~\ref{thm:sparse} are in \ref{app:sparse}: a qualifying denominator has a continued fraction with total excess $\sum_i(a_i-1)\le\frac1{10}\log p+O(1)$, and such words are too few.

Both constants are far from optimal. Applying the three-point inequality at every move, including those away from the steps with negative residual, gives a much larger penalty.

\begin{theorem}[Computer-assisted]\label{thm:penalty}
For every unit $r\le p/2$ modulo $p\ge2$, with $p/r=[a_1;\dots,a_K]$ canonical and $E=\sum_i(a_i-1)$,
\[
\frac{\Phi(p,r)}{p^2}\ge C_{\mathrm g}\log p+\frac{E-1}{250}-D_0+\frac1{4p},\qquad D_0=\frac45C_{\mathrm g}\log2-\frac1{40}<\frac{277}{1000}.
\]
Consequently, as $X\to\infty$,
\begin{align*}
\#\{2\le p\le X:\ m(p)\le C_{\mathrm g}\,p^2\log p\}&=O\bigl((\log X)^{69}\bigr),\\
\#\{2\le p\le X:\ m(p)\le(C_{\mathrm g}+\tfrac1{200000})\,p^2\log p\}&=O(X^{1/90}),\\
\#\{2\le p\le X:\ m(p)\le(C_{\mathrm g}+\tfrac1{750})\,p^2\log p\}&=O(X^{39/40}).
\end{align*}
In particular $m(p)>(C_{\mathrm g}+\frac1{750})p^2\log p$ for a set of integers $p$ of density one and for a set of primes of relative density one.
\end{theorem}

The proof (\ref{app:penalty}) bounds the three-point sums $\varrho_{j-1}+3\varrho_j+\varrho_{j+1}$ below by $0$, $1/40$, or $1/100$ according as the incoming and outgoing moves are two flips, two no-flips, or one of each. The last case is the computer-assisted part: it reduces to the positivity of one explicit polynomial on a rectangle, certified by exact rational Bernstein coefficients.

\begin{remark}[Scale]\label{rem:scale}
The counts of Theorems~\ref{thm:sparse} and~\ref{thm:penalty} are asymptotic and are informative only for very large $X$. The final majorant $Le^{1/31}X^\theta(32e)^{D_1}$ of the count behind Theorem~\ref{thm:sparse} in \ref{app:sparse} exceeds $X$ for every $2\le X\le10^{900000}$, and the explicit bound behind the first count of Theorem~\ref{thm:penalty} in \ref{app:penalty}, $70+\sum_{2\le k\le\mathsf A\log X}\binom{k+68}{68}$ with $\mathsf A=1/\log\tau$, exceeds $X$ for every $2\le X\le10^{83}$. At accessible sizes the typical behavior is not yet visible: among $2\le p\le3000$ only $p=506$ has $m(p)>C_{\mathrm g}p^2\log p$, and even there the ratio, $0.54531\ldots$, is below the lower end $C_{\mathrm g}+\frac1{750}=0.54538\ldots$ of Theorem~\ref{thm:typical-upper}; moreover $m(p)/(p^2\log p)=0.5405\ldots$ at $p=221688$ (Proposition~\ref{prop:counterexample}), below $C_{\mathrm g}$. The set bounded by the first count of Theorem~\ref{thm:penalty} is infinite and has at least $(1/\log\tau+o(1))\log X$ elements up to $X$: since $c_{\mathrm F}<0$ in Remark~\ref{rem:fib-second-order}, $m(F_n)\le\Phi(F_n,F_{n-1})<C_{\mathrm g}F_n^2\log F_n$ for all large $n$. Continued fractions with further isolated digits $2$ appear to qualify as well. For $p/r=[2;1^{40},2,1^{40},2]$ one finds $\Phi(p,r)/p^2-C_{\mathrm g}\log p=-0.22044\ldots$; each further digit $2$, thirty digits $1$ from the others, raises this value by $0.02188\ldots$, and with eleven interior digits $2$ it is still $-0.00163\ldots$. If such words give distinct values of $p$ up to reversal, this suggests at least order $(\log X)^{12}$ elements. The true order is Problem~\ref{prob:exceptional} of Section~\ref{sec:open}.
\end{remark}

\subsection{The typical size of the minimum}

\begin{remark}[An upper bound for every $p$]\label{rem:abh-upper}
Larcher \cite[Cor.~2]{Larcher1986} proved $\Sigma_{\min}(p)=O\bigl((p/|\Units p|)\log p\log\log p\bigr)$. A tail estimate of Rukavishnikova \cite{Rukavishnikova2006}, quoted in the introduction of \cite{Rukavishnikova2011} and \cite[Sect.~1, (2)]{AistleitnerBordaHauke2024}, removes the factor $p/|\Units p|$, and Aistleitner, Borda, and Hauke \cite[Cor.~2]{AistleitnerBordaHauke2024} prove that every $N\ge3$ has a unit $a$ with canonical digit sum $\Sigma(N,a)\le\frac{12}{\pi^2}\log N\log\log N+O(\log N)$. The upper bound of Theorem~\ref{thm:pq} and \eqref{eq:length}, applied to this unit, give
\[
m(p)\le\frac{2}{\pi^2}\,p^2\log p\log\log p+O(p^2\log p)\qquad(p\ge3).
\]
Theorem~\ref{thm:typical-upper} below is much stronger on a set of density one, but this is the strongest unconditional upper bound we know that holds for every $p$.
\end{remark}

On the upper side, the density-one method of Bourgain and Kontorovich \cite{BourgainKontorovich2014} applies: Huang's theorem \cite[Thm.~1.6]{Huang2015}, with \cite[Thm.~4.12]{PollicottVytnova2022}, gives for almost all $p$ a unit with all digits at most $5$, and Theorem~\ref{thm:pq} then gives a constant of about $2.605$. The next theorem lowers it to $0.584$.

\begin{theorem}[Computer-assisted]\label{thm:typical-upper}
There is $c_0>0$ such that
\[
\#\Bigl\{2\le p\le X:\ m(p)\ge\frac{73}{125}\,p^2\log p\Bigr\}=O\bigl(Xe^{-c_0\sqrt{\log X}}+X^{1-1/1250}\bigr).
\]
Consequently, for all $p\le X$ outside a set of this size, and so for a set of integers of density one and a set of primes of relative density one,
\[
C_{\mathrm g}+\frac1{750}<\frac{m(p)}{p^2\log p}<\frac{73}{125}=0.584 .
\]
\end{theorem}

The proof (\ref{app:typical}) has three ingredients. First, an exact formula expresses $\Phi(p,r)/p^2$ as a sum over the digits of any continued fraction of $r/p$. Second, a transfer-matrix certificate on words of eight digits from $\{1,2,3,4\}$ shows that the pairs $(p,r)$ with such digits, $p\le X$, and $\Phi(p,r)\ge\frac{73}{125}p^2\log p$ number $O(X^{197/125})$; this is the computer-assisted part. Third, two theorems are imported: Kan's theorem \cite[Thm.~1.4, formulas (1.6)--(1.8)]{Kan2021} that almost all admissible denominators $p$ have $\gg p^{2\delta-1.001}$ numerators with digits in an alphabet whose continued fraction Cantor set has Hausdorff dimension $\delta\ge(\sqrt{40}-4)/3$, and the rigorous bound $\delta>0.7889$ of Pollicott and Vytnova \cite[Thm.~4.14]{PollicottVytnova2022} for the alphabet $\{1,2,3,4\}$, for which every integer is admissible because it is for $\{1,2\}$.

\begin{remark}[Conditional finiteness]\label{rem:zhang}
If there is an absolute $\mathcal M$ such that every $p\ge2$ has a unit $r$ whose continued fraction (canonical or ending in $1$) has all digits at most $\mathcal M$, then $\Sigma(p,r)\le\mathcal M(K+1)$, and Theorem~\ref{thm:pq} with \eqref{eq:length} gives
\[
\limsup_{p\to\infty}\frac{m(p)}{p^2\log p}\le\frac{\mathcal M/6+121/288}{\log\tau}<\infty .
\]
This hypothesis is Zaremba's conjecture. A 2026 preprint of Shkredov \cite[Cor.~1]{Shkredov2026} states it for prime and square-free denominators, and for all denominators in its second version \cite[Remark~9]{Shkredov2026}; a 2026 preprint of Zhang \cite[Thm.~1.2]{Zhang2026}, building on it, states it for all denominators, with a large absolute $\mathcal M$. The extension in \cite[Remark~9]{Shkredov2026} rests on work of Zhang and Tang, so the two claims are not independent. We have verified neither proof. Along the primes and along square-free $p$ the finiteness of the limsup rests on \cite[Cor.~1]{Shkredov2026} alone; along all $p$ it rests on both preprints. Even then the value of the limsup is not determined; it is at least $C_{\mathrm g}+\frac1{750}$ by Theorem~\ref{thm:penalty}.
\end{remark}

\section{The minimizers}\label{sec:minimizers}

The next statement is elementary and reduces the finiteness of $\limsup m(p)/(p^2\log p)$ to the minimal digit sum; Remark~\ref{rem:abh-upper} gives the known upper bounds for that sum.

\begin{proposition}\label{prop:order}
For $p\ge2$ and every unit $r$ modulo $p$,
\[
\frac{p^2}8\bigl(\Sigma(p,r)-1\bigr)+\frac p4\le\Phi(p,r)\le p^2\Bigl(\frac{\Sigma(p,r)}2-\frac34\Bigr)-\frac p4+\frac12 .
\]
Consequently $\limsup_{p\to\infty}m(p)/(p^2\log p)<\infty$ if and only if $\Sigma_{\min}(p)=O(\log p)$, a conjecture attributed to Moser (see \cite[Sect.~1]{Shkredov2026}), and $m(p)=p^2\Sigma_{\min}(p)/6+O(p^2\log p)$.
\end{proposition}

\begin{proof}
For $x,y\ge1$, $\frac18\le\Gamma(x,y)\le\frac12$: the lower bound was shown in \ref{app:sparse}, and with $a=x-1$, $b=y-1$ one has $2(x+y-1)^2-x^2-y^2=a^2+b^2+4ab+2a+2b\ge0$. Insert these bounds, the step count $\Sigma(p,r)-2$, and the bounds on the boundary term of Proposition~\ref{prop:action}. Minimizing over $r$ gives the equivalence, and Theorem~\ref{thm:pq} with \eqref{eq:length} gives the last formula.
\end{proof}

The minimizers are constrained: the next proposition, proved in \ref{app:minimizers}, bounds every digit, the digit sum, and the residue of a minimizer in terms of $m(p)/p^2$.

\begin{proposition}\label{prop:minimizers}
Let $r$ be a unit modulo $p\ge2$ with inverse $r^*\in[1,p-1]$, and let $H\ge\Phi(p,r)$ be a real number. Then $\Sigma(p,r)\le1+\lfloor(8H-2p)/p^2\rfloor$, every digit satisfies $a_i<6H/p^2+4$, the first and last digits satisfy $a_1,a_K<6H/p^2+2$, and
\[
\min\{r,\ p-r,\ r^*,\ p-r^*\}>\frac{p}{6H/p^2+3}.
\]
\end{proposition}

The minimizers are not determined by the digit sum, however, and need not minimize it.

\begin{proposition}[Computer-assisted]\label{prop:counterexample}
For $p=221688$, $m(p)=326994450836$, attained exactly at $r\in\{65003,91645,130043,156685\}$, each of digit sum $29$; and $\Sigma_{\min}(p)=28$, attained exactly at $r\in\{84859,84883,136805,136829\}$, each with $\Phi(p,r)=327143115683$. In particular the two sets of minimizers are disjoint.
\end{proposition}

The proof is an exhaustive evaluation of $\Phi(p,r)$ over all $73872$ units by Theorem~\ref{thm:closed-form}, confirmed by a second implementation of a different exact formula and by a direct evaluation of the two representative values from Definition~\ref{def:phi}; the computations are described in Section~\ref{sec:computation}. For example $65003/221688=[0;3,2,2,3,2,3,2,2,2,1,1,2,2,2]$ and $84859/221688=[0;2,1,1,1,1,2,1,1,1,1,1,1,1,4,1,2,1,1,1,1,2]$.

\section{Open problems}\label{sec:open}

\begin{enumerate}
\item \emph{The upper behavior of the minimum.} Determine the limit superior of $m(p)/(p^2\log p)$. The ratio is $O(\log\log p)$ (Remark~\ref{rem:abh-upper}), and the limit superior is finite if the two linked preprints discussed in Remark~\ref{rem:zhang} are correct; its value would remain open. Is there a constant $A$ with $m(p)\sim Ap^2\log p$ for a set of integers $p$ of density one, or for a set of primes of relative density one, and is it the same constant? Theorems~\ref{thm:penalty} and~\ref{thm:typical-upper} place such an $A$ in $[C_{\mathrm g}+1/750,\ 0.584]$.
\item \emph{The liminf along the primes.} Is $\liminf m(\ell)/(\ell^2\log\ell)=C_{\mathrm g}$ as $\ell$ runs through the primes? Theorem~\ref{thm:golden} attains $C_{\mathrm g}$ along the Fibonacci numbers, and Theorems~\ref{thm:penalty} and~\ref{thm:typical-upper} concern typical primes.
\item \emph{The Fibonacci minimizers.} Is $m(F_n)=\Phi(F_n,F_{n-1})$ for all $n\ge3$? This holds for $3\le n\le28$. More generally, describe the residues $r$ attaining $m(p)$; minimizing the sum of partial quotients does not select them (Proposition~\ref{prop:counterexample}).
\item\label{prob:exceptional} \emph{The exceptional set.} Determine the order of $\#\{2\le p\le X:\ m(p)\le C_{\mathrm g}\,p^2\log p\}$. It is at least of order $\log X$ (Remark~\ref{rem:scale}) and $O((\log X)^{69})$ (Theorem~\ref{thm:penalty}), and the computations in Remark~\ref{rem:scale} suggest at least order $(\log X)^{12}$. Is it $(\log X)^{d+o(1)}$ for some $d$, and what is $d$?
\item \emph{The second-order term.} Put $D(p)=m(p)/p^2-C_{\mathrm g}\log p$. By Theorem~\ref{thm:golden} and Remark~\ref{rem:fib-second-order},
\[
\tfrac18-C_{\mathrm g}\log2\le\liminf_{n\to\infty}D(F_n)\le\limsup_{n\to\infty}D(F_n)\le c_{\mathrm F}<0 ;
\]
the limit $c_{\mathrm F}$ is proved for $\Phi(F_n,F_{n-1})$, and for $m(F_n)$ it would follow from $\Phi(F_n,F_{n-1})-m(F_n)=o(F_n^2)$. Is $\liminf_{p\to\infty}D(p)=c_{\mathrm F}$? On a set of density one, $D(p)\to\infty$ by Theorem~\ref{thm:sparse}.
\end{enumerate}

\section{Computation and reproducibility}\label{sec:computation}

The code \cite{code2026} reproduces each computational claim of the paper. It was written with the assistance of Claude Opus (Anthropic), as stated in the declaration at the end of the paper. It uses only the Python standard library; the verification uses exact integers and rationals and rational interval bounds for logarithms (only the optional script that carries out the search for the vector of \ref{app:typical} uses floating point). The file \texttt{COMPUTATIONS.md} of the code lists, for each statement, the script that checks it and the range of the finite tests; these accompany the proofs and replace none of them. The outputs are deterministic, and a rerun is compared byte for byte with the shipped manifest of SHA-256 digests. The checks run in parallel; the full suite takes about a minute on two cores and gives identical outputs under Python 3.11, 3.12, 3.13, and 3.15, and the quick suite also under Python 3.9.

Three results are computer-assisted, and the code includes their certificates and verifiers: Proposition~\ref{prop:counterexample}, by an exact evaluation of $\Phi(221688,r)$ for all $73872$ units with two different formulas; Theorem~\ref{thm:penalty}, by $1296$ Bernstein coefficients of one polynomial, all positive, computed by two implementations; and Theorem~\ref{thm:typical-upper}, by a positive integer vector for a transfer matrix with $65536$ states and $262144$ edges, verified on all rows and, by a separately written program, on $163$ rows. The last also relies on the theorems of Kan \cite{Kan2021} and of Pollicott and Vytnova \cite{PollicottVytnova2022}, whose proofs the code does not check. The digits of $c_{\mathrm F}$ in Remark~\ref{rem:fib-second-order}, the range in Remark~\ref{rem:pq-range}, the equality $m(F_n)=\Phi(F_n,F_{n-1})$ for $3\le n\le28$, and the values in Remark~\ref{rem:scale} are interval or exhaustive evaluations, used in no proof; the quick suite runs the same checks on smaller ranges (for instance $n\le18$ for the Fibonacci minimizers). All other results have complete proofs in the text; the code also checks the polynomial and rational identities used in them, symbolically where it works with polynomials and otherwise at sample rational points.

\appendix
\renewcommand{\thetheorem}{\Alph{section}.\arabic{theorem}}
\section{Proofs for Section~\ref{sec:reciprocity}}\label{app:two}

\subsection{Proof of Theorem~\ref{thm:closed-form}}\label{app:closed-form}

We prove (i), with the convention $Q_{1,0}\equiv0$, for $q=cp+s$, $0\le s<p$:
\[
Q_{q,p}(x,y)=P(p,s,c)\,x^2+Q_{p,s}(cx+y,\,x).
\]
If $s>0$, expand the right side: its coefficients of $x^2$, $xy$, $y^2$ are $P(p,s,c)+c^2\Phi(p,s)+2c\M(s,p)+\Phi(s,p)$, $2c\Phi(p,s)+2\M(s,p)$, and $\Phi(p,s)$. They equal $\Phi(q,p)$, $2\M(p,q)$, and $\Phi(p,q)$ by Lemma~\ref{lem:integrated}, Theorem~\ref{thm:euclid-step}, and periodicity. If $s=0$ then $p=1$, $\M(1,q)=\Phi(1,q)=0$, and \eqref{eq:form-step-main} reads $\Phi(q,1)x^2=P(1,0,q)x^2$, which is the terminal case of Lemma~\ref{lem:integrated}. This is the last Euclidean step, $i=K$.

At step $i$ the arguments of the remaining form $Q_{n_{i-1},n_i}$ are $X_i=u_ix+v_iy$ and $Y_i=u_{i-1}x+v_{i-1}y$: this holds for $i=1$, and by \eqref{eq:form-step-main} the next arguments are $c_iX_i+Y_i=X_{i+1}$ and $X_i=Y_{i+1}$. Iterating \eqref{eq:form-step-main} for $i=1,\dots,K$ gives
\begin{equation}
Q_{q,p}(x,y)=\sum_{i=1}^KP_i\,(u_ix+v_iy)^2,
\end{equation}
which is (ii); comparing coefficients gives the three sums.
\qed

\section{Proofs for Sections~\ref{sec:pq} to~\ref{sec:minimizers}}\label{app:extremes}

We keep the notation of Sections~\ref{sec:prelim} and~\ref{sec:golden}. For a unit $r$ modulo $p\ge2$ let $(n_i)$, $(c_i)=(a_i)$, and $(u_i)$ be the data of Theorem~\ref{thm:closed-form} for the chain $(n_0,n_1)=(p,r)$, and let $q_{-1}=0$, $q_0=1$, $q_i=a_iq_{i-1}+q_{i-2}$ be the continuants of $r/p$, so $u_i=q_{i-1}$ and $q_K=p$. Induction on $i$ gives the determinant identity
\begin{equation}\label{eq:det-identity}
n_{i-1}q_{i-1}+n_iq_{i-2}=p\qquad(1\le i\le K+1),
\end{equation}
and the classical congruences $q_{i-1}r\equiv(-1)^{i-1}n_i\pmod p$ for $0\le i\le K+1$.

\subsection{Proof of Theorem~\ref{thm:pq}}\label{app:pq}

\emph{Upper bound and the refined lower bound.} By Theorem~\ref{thm:closed-form}, $\Phi(p,r)=\sum_iu_i^2P(n_i,n_{i+1},a_i)$. Fix $i$ and write $a=a_i$, $n=n_i$, $s=n_{i+1}<n$, $u=u_i$, $v=u_{i-1}\le u$, and
\[
X=\frac{un}p,\qquad Y=\frac{us}p,\qquad Z=\frac{vn}p .
\]
By \eqref{eq:det-identity}, $aX+Y+Z=1$, and $0\le Y,Z\le X$. Expanding the closed form of $P$ and separating the terms of degree two in $(n,s)$ from those of degree one gives
\[
\frac{u^2P(n,s,a)}{p^2}=\frac a6(aX+Y+Z)^2+\mathcal E+\mathcal L=\frac a6+\mathcal E+\mathcal L,
\]
\begin{gather*}
\mathcal E=\frac a{12}\bigl(3aX^2+2aXY-4aXZ+X^2+6XY+4Y^2\\-4YZ-2Z^2\bigr),\\
\mathcal L=-\frac1n\Bigl(\frac{a(a+1)X^2}4+\frac{aXY}2\Bigr).
\end{gather*}
Put $\eta=Y/X$, $\zeta=Z/X\in[0,1]$, and $\Delta=a+\eta+\zeta$, so $X=1/\Delta$. Then $\mathcal E=a\mathcal N/(12\Delta^2)$ with $\mathcal N=3a+1+(2a+6)\eta-4a\zeta+4\eta^2-4\eta\zeta-2\zeta^2$, and
\begin{gather*}
\Delta\mathcal N_\eta-2\mathcal N=2(a^2-1)+6(a-1)\eta+6(a+1)\zeta+12\eta\zeta\ge0,\\
\Delta\mathcal N_\zeta-2\mathcal N=-\bigl(4a^2+6a+2+12(a+1)\eta+12\eta^2\bigr)<0 .
\end{gather*}
Since $\partial_\eta(\mathcal N/\Delta^2)=(\Delta\mathcal N_\eta-2\mathcal N)/\Delta^3$ and similarly for $\zeta$, $\mathcal E$ is nondecreasing in $\eta$ and decreasing in $\zeta$, so its extreme values on $[0,1]^2$ are at $(\eta,\zeta)=(0,1)$ and $(1,0)$:
\[
-\frac a{12(a+1)}\le\mathcal E\le\frac{a(5a+11)}{12(a+1)^2}\le\frac{121}{288},
\]
the last because $\frac{d}{da}\frac{5a^2+11a}{(a+1)^2}=\frac{11-a}{(a+1)^3}$, so the maximum over real $a\ge1$ is at $a=11$. Next $\mathcal L\le0$, and since $Y\le1-aX$ and $0\le X\le1/a$, the bracket in $\mathcal L$ is at most $f(X)=a(1-a)X^2/4+aX/2$, which is increasing on $[0,1/a]$ (as $f'(X)=\frac a2(1-(a-1)X)\ge\frac12$); hence $0\le-\mathcal L\le f(1/a)/n=(a+1)/(4an)\le1/(2n)$. The remainders satisfy $n_K=1$, $n_{K-1}\ge2$, and $n_{i-1}\ge n_i+n_{i+1}$, so $n_{K-j}\ge F_{j+2}$, and as $F_{j+2}\ge2F_j$,
\[
\begin{split}\sum_{i=1}^K\frac1{n_i}\le\sum_{j\ge2}\frac1{F_j}&\le\sum_{j=2}^{20}\frac1{F_j}+\frac2{F_{21}}+\frac2{F_{22}}\\&=\frac{21724319459128235137204000223387}{9205445844772861394529630892080}<\frac{59}{25}.\end{split}
\]
Summing over $i$ gives $-p^2(K/12+59/50)\le\Phi(p,r)-p^2\Sigma(p,r)/6\le\frac{121}{288}Kp^2$.

\emph{The lower bound $-Kp^2/6$.} For $1\le k<p$ let $d_k^-$ and $d_k^+$ be as in Section~\ref{sec:spacings} for $(n,a)=(p,r)$; by \eqref{eq:insertion}, $\Phi(p,r)=\sum_{k=1}^{p-1}k\,d_k^-d_k^+$.

\emph{Records.} Suppose $U,V,A,B\ge1$, $Ur\equiv A$, $Vr\equiv-B\pmod p$, and $AV+BU=p$. Then $d_k^-=A$ and $d_k^+=B$ for $\max(U,V)\le k<U+V$. Indeed $(U,A)$ and $(V,-B)$ lie in the lattice $\{(j,t)\in\Z^2:t\equiv jr\pmod p\}$ of index $p$, and their determinant is $-p$, so they form a basis. A lattice point $(j,t)=m(U,A)+m'(V,-B)$ has $j\ge U+V$ if $m,m'\ge1$; $j\le0$ if $m,m'\le0$; $t\ge A$ if $m\ge1$ and $m'\le0$; and $t\le-B$ if $m\le0$ and $m'\ge1$. Hence no $j$ with $0<j<U+V$ has a representative $t\equiv jr$ with $-B<t<A$, while the indices $U,V\le k$ realize $A$ and $B$.

For $1\le i\le K$ and $1\le j\le a_i$ the congruences above give $q_{i-1}r\equiv(-1)^{i-1}n_i$ and $(q_{i-2}+jq_{i-1})r\equiv(-1)^i(n_{i-1}-jn_i)$, and \eqref{eq:det-identity} gives $(n_{i-1}-jn_i)q_{i-1}+n_i(q_{i-2}+jq_{i-1})=p$. Applying the record statement (with the roles of the two lengths exchanged when $i$ is odd), $d_k^-d_k^+=n_i(n_{i-1}-jn_i)$ for $q_{i-2}+jq_{i-1}\le k<q_{i-2}+(j+1)q_{i-1}$, whenever $n_{i-1}-jn_i\ge1$. These ranges tile $[1,p+q_{K-1})$ consecutively, block $i$ covering $[q_{i-2}+q_{i-1},q_i+q_{i-1})$; only $[1,p)$ is needed, and the only range with $n_{i-1}-jn_i=0$ is $(i,j)=(K,a_K)$, which starts at $q_K=p$. For example, $5/13=[0;2,1,1,2]$ has remainders $(n_0,\dots,n_4)=(13,5,3,2,1)$ and continuants $(q_0,\dots,q_4)=(1,2,3,5,13)$; the four blocks cover $[1,3)$, $[3,5)$, $[5,8)$, and $[8,18)$, and in the first block $d_1^-d_1^+=5\cdot8$ and $d_2^-d_2^+=5\cdot3$, in agreement with $\phi_{13,5}(1)-\phi_{13,5}(2)=40$ and $\phi_{13,5}(2)-\phi_{13,5}(3)=15$. Therefore $\Phi(p,r)=\sum_{i=1}^KC_i$ with
\[
C_i=\frac{n_iq_{i-1}}2\sum_{j=1}^{a_i}(n_{i-1}-jn_i)\bigl(2q_{i-2}+(2j+1)q_{i-1}-1\bigr).
\]
Fix $i$, write $a=a_i$, and put $x=q_{i-2}/q_{i-1}$, $y=n_{i+1}/n_i$, $\delta=1/q_{i-1}$, $z=n_iq_{i-1}$; then $x,y\in[0,1]$, $\delta\in(0,1]$, and $p=z(a+x+y)$ by \eqref{eq:det-identity}. Summing the powers of $j$,
\[
\frac{C_i}{z^2}=\frac a{12}\Bigl(6(a+y)(2x+a+2-\delta)-(a+1)(6x+4a+5-3\delta)\Bigr),
\]
so $C_i-ap^2/6=z^2\mathcal D$ with
\[
\mathcal D=a^2\Bigl(\frac{x+y}6+\frac{1-\delta}4\Bigr)+a\Bigl(\frac{2xy}3-\frac{x^2+y^2}6+\Bigl(1-\frac\delta2\Bigr)y-\frac x2-\frac{5-3\delta}{12}\Bigr).
\]
Since $\partial_\delta\mathcal D=-a(a-1)/4-ay/2\le0$, it suffices to bound $\mathcal D$ below at $\delta=1$ (the upper bound of Theorem~\ref{thm:pq} was proved above). At $\delta=1$, $G(a)=6\mathcal D+(a+x+y)^2=a^2(1+x+y)+a(4xy-x^2-y^2+5y-x-1)+(x+y)^2$ has $G(1)=6y(1+x)\ge0$, $G'(1)=1+x+7y+4xy-x^2-y^2\ge0$, and $G''>0$, so $G(a)\ge0$ for $a\ge1$. Hence $C_i-a_ip^2/6\ge-p^2/6$, and summing over the $K$ blocks gives $\Phi(p,r)-p^2\Sigma(p,r)/6\ge-Kp^2/6$. The final assertion of Theorem~\ref{thm:pq} follows from \eqref{eq:length}. \qed

\subsection{Proof of Proposition~\ref{prop:action}}\label{app:action}

\emph{Spacings along the path.} Consider a step $(A,B,U,V)$, $A>B$. By induction along the algorithm, $\gcd(U,V)=1$, $AV+BU=p$, and, for a sign $\epsilon\in\{\pm1\}$ that changes at each interchange, $U\epsilon r\equiv A$ and $V\epsilon r\equiv-B\pmod p$; initially $(U,V)=(1,1)$, and the update $(A,U)\mapsto(A-B,U+V)$ preserves all three. Replacing $r$ by $\epsilon r$ does not change any spacing (Proposition~\ref{prop:sym}), so take $\epsilon=1$. Stage $h=U+V$ consists of the points $ir$, $0\le i<U+V$. The map $i\mapsto i+U\bmod (U+V)$ is a single cycle, since $\gcd(U,V)=1$; it moves $ir$ forward by $Ur\equiv A$ if $i<V$ and by $-Vr\equiv B$ if $i\ge V$. These $U+V$ positive displacements add up to $VA+UB=p$, so the cycle winds once around $\Z/p\Z$ and lists the points in circular order: stage $U+V$ has $V$ spacings of length $A$ and $U$ of length $B$. For $0\le t<V$ the point $(U+V+t)r$ lies at displacement $A-B\in(0,A)$ after the point $tr$, inside the spacing of length $A$ that starts at $tr$; these $V$ spacings are distinct. So the next $V$ insertions each split a spacing $A$ into $A-B$ and $B$, lowering $\phi$ by $B(A-B)$, and lead to the next step. Summing $\phi$ over the stages $h=U+V,\dots,U+2V-1$ gives
\[
\begin{split}Q&=V\Bigl(V\binom A2+U\binom B2\Bigr)-B(A-B)\binom V2\\&=\frac12\Bigl(V^2(A-B)^2+pVB-pV+VB(A-B)\Bigr).\end{split}
\]

\emph{Summation.} The steps cover the stages $h=2,\dots,p-1$, the step $j$ covering $V_j$ of them, so $\sum_jV_j=p-2$ and $\sum_jV_jB_j(A_j-B_j)=\phi_{p,r}(2)=\binom r2+\binom{p-r}2$. Put $t=AV/p$ and $z=BV/p$. The next step has $t'=t-z$, or $t'=1-(t-z)$ after an interchange (using $AV+BU=p$), and in both cases $(t-z)^2+z=t^2+t(1-t)-t'(1-t')$. The first step has $t(1-t)=r(p-r)/p^2$ and the terminal state $(1,1,U_*,V_*)$ has $t(1-t)=U_*V_*/p^2$. Hence $\sum_j\frac12(V_j^2(A_j-B_j)^2+pV_jB_j)=\frac12\sum_j(A_jV_j)^2+\frac12(r(p-r)-U_*V_*)$, and with $\Phi(p,r)=\binom p2+\sum_jQ_j$,
\begin{equation}\label{eq:action-one}
\Phi(p,r)=\frac{p^2+p}4-\frac{U_*V_*}2+\frac12\sum_j(A_jV_j)^2 .
\end{equation}

\emph{The reversed path.} To a step $s=(A,B,U,V)$ associate $s'=(U+V,V,A-B,B)$, which again satisfies $A'V'+B'U'=p$ and $A'>B'$. If $s_{j+1}$ follows $s_j$ without interchange, $s_{j+1}'=(U+2V,V,A-2B,B)$, and one subtraction sends it to $s_j'$; if with interchange, $s_{j+1}'=(U+2V,U+V,2B-A,A-B)$, and one subtraction followed by an interchange sends it to $s_j'$. The last step of the original path has lengths $(2,1)$, so its image has lengths $(U_*,V_*)$ and counts $(1,1)$, and the image of the first step $(A,B,1,1)$ is followed by the terminal state with counts $(A,B)=\{r,p-r\}$. Hence, if $s_1,\dots,s_N$ are the steps of $r$, then $s_N',\dots,s_1'$ is the path of a unit $\tilde r\in\{U_*,V_*\}$, and the terminal congruences $U_*r\equiv\pm1$ give $\tilde r\equiv\pm r^{-1}$, so $\Phi(p,\tilde r)=\Phi(p,r)$ by Proposition~\ref{prop:sym}. Formula \eqref{eq:action-one} for this path reads $\Phi(p,r)=\frac{p^2+p}4-\frac{r(p-r)}2+\frac12\sum_j((U_j+V_j)B_j)^2$. Average the two formulas. Since $p=AV+BU=BV(x+y-1)$, we have $AV/p=x/(x+y-1)$ and $(U+V)B/p=y/(x+y-1)$, which gives the first identity of Proposition~\ref{prop:action}. For $p=2$ the path is empty and the identity reads $1/4=1/4$.

\emph{Logarithms, boundary, and step counts.} A step changes the length sum from $A+B$ to $A$ and the count sum from $U+V$ to $U+2V$; these ratios are $1+1/x$ and $1+1/y$, and the length sums run from $p$ to $2$ while the count sums run from $2$ to $p$. This proves $\sum_j\Lambda=\log(p/2)$. Both $U_*V_*$ and $r(p-r)$ are products of positive integers with sum $p$, so they lie in $[p-1,p^2/4]$, which gives the bounds on the boundary term. Let $r\le p/2$ and $p/r=[a_1;\dots,a_K]$. The ordinary subtractive algorithm on $(p,r)$ performs $\Sigma(p,r)$ subtractions; the path omits the first one, $(p,r)\mapsto(p-r,r)$, and the last one, $(1,1)\mapsto(1,0)$, so it has $\Sigma(p,r)-2$ steps. If $K=1$, then $r=1$ and the $p-2$ steps $(p-1,1),\dots,(2,1)$ are no-flips. If $K\ge2$, a step inside block $i<K$ has $x=a_i-j+n_{i+1}/n_i$ for $0\le j<a_i$ ($1\le j<a_1$ for $i=1$), which is less than $2$ only for $j=a_i-1$, while the $a_K-1$ steps of the last block have $x\ge2$; so there are $K-1$ flips and $\Sigma(p,r)-K-1=E-1$ no-flips. \qed

\subsection{Proofs of Lemmas~\ref{lem:outside} and~\ref{lem:stencil}}\label{app:golden}

\emph{Constants.} From $2.236067^2<5<2.236068^2$ and $\tau>(1+2.236067)/2$, the series $\log\frac{1+a}{1-a}=2\sum_{j\ge0}a^{2j+1}/(2j+1)$ with $a=1236067/5236067$ gives $\log\tau>2(a+a^3/3+a^5/5+a^7/7)>0.4812$, an exact rational comparison, and hence $C_{\mathrm g}<5.236068/(20\cdot0.4812)<\kappa=5441/10000$. Also $\log\tau>\log\frac32=2\operatorname{artanh}\frac15>\frac25$ and $3+\sqrt5<6$ give $C_{\mathrm g}<\frac34<\frac45$.

\emph{A rational majorant.} For $u\ge1$ put $w=1/(2u+1)\in(0,1/3]$. Then $\log(1+1/u)=2\sum_{j\ge0}w^{2j+1}/(2j+1)<2w+\frac{2w^3}{3(1-w^2)}=f(u)$, with
\[
f(u)=\frac16\Bigl(\frac1u+\frac1{u+1}+\frac8{2u+1}\Bigr)=\frac{12u(u+1)+1}{6u(u+1)(2u+1)} .
\]

\begin{proof}[Proof of Lemma~\ref{lem:outside}]
It suffices to show $\Gamma(x,y)\ge\kappa(f(x)+f(y))/2$ for $x,y\ge1$ with $\max(x,y)\ge2$. Put $s=x+y\ge3$ and $t=xy$. A direct computation gives
\[
\begin{split}\frac{\Gamma(x,y)}{(f(x)+f(y))/2}&=\frac{3(s^2-2t)\,t\,(t+s+1)(4t+2s+1)}{(s-1)^2(s+1)\bigl(24t^2+6(4s+3)t+s(2s+1)\bigr)}\\&=:\mathcal R_s(t).\end{split}
\]
Writing $\mathcal R_s=3\mathcal N/((s-1)^2(s+1)\mathcal D)$, the sign of $\mathcal R_s'(t)$ is that of $\mathcal N'\mathcal D-\mathcal N\mathcal D'$, a quintic in $t$ with coefficients, from degree five down to zero,
\begin{gather*}
-384,\quad 96s^2-864s-672,\quad 192s^3-496s^2-944s-360,\\
120s^4-246s^2-186s-36,\quad 24s^5+16s^4-22s^3-20s^2-4s,\\
4s^6+8s^5+5s^4+s^3 .
\end{gather*}
For $s\ge2$ the last three are positive (for degree two use $20s^4\ge80s^2$, so one sixth of it is at least $39s^2-31s-6$; for degree one, half of it divided by $s$ is at least $8s^3+37s^2-10s-2$). If the degree-four coefficient is positive then $s>9$, and then the degree-three coefficient is $8s^2(24s-62-118/s-45/s^2)>0$. So the coefficients have exactly one sign change, the quintic has exactly one positive root, and $\mathcal R_s$ increases and then decreases in $t>0$; its minimum over an interval is at an endpoint. For $3\le s\le4$ the admissible products form $[s-1,2(s-2)]$, and for $s\ge4$ they form $[s-1,s^2/4]$. Hence it suffices to treat three boundaries.

On $x=y=z\ge2$, clearing positive denominators turns $\Gamma-\kappa(f(x)+f(y))/2\ge0$, up to a positive constant factor, into $P_d(z)=60000z^5-171168z^4+30000z^3+174112z^2-43528z-5441\ge0$, and
\[
\begin{split}P_d(2+v)&=25263+335544v+1046080v^2+1060656v^3\\&\quad+428832v^4+60000v^5>0\qquad(v\ge0).\end{split}
\]
On $x=1$, $y=z\ge2$, multiplying by $720000z^2(z+1)(2z+1)$ gives $P_1(z)=87950z^4-259827z^3+12223z^2+507354z+180000$, and $P_1(2+v)=572184+252722v+564061v^2+443773v^3+87950v^4>0$. On $x=2$, $y=z\in[1,2]$, with $f(2)=73/180$, multiplying by $3600000z(z+1)^2(2z+1)$ gives
\[
P_2(z)=1005614z^4-3044725z^3+1693708z^2+1080817z-163230 .
\]
With $v=z-5/3\in[-2/3,1/3]$, $P_2=\frac{509840}{81}-\frac{640346}{27}v+\frac{9690949}3v^2+\frac{10978105}3v^3+1005614v^4$. For $v\ge0$ drop the last two terms; the remaining quadratic has discriminant $-58879958258204/729<0$. For $v=-w\le0$, $P_2=\frac{509840}{81}+\frac{640346}{27}w+w^2\bigl(\frac{9690949}3-\frac{10978105}3w+1005614w^2\bigr)$, and the parenthesis is at least $\frac{9690949}3-\frac23\cdot\frac{10978105}3=\frac{7116637}9>0$ for $w\le2/3$. Thus $\Gamma>\kappa\Lambda$ on all three boundaries, hence everywhere in the region.
\end{proof}

\begin{proof}[Proof of Lemma~\ref{lem:stencil}]
Let $\mathcal S(x,y)$ be the left side and $s=x+y$, $t=xy$. The three rational terms have the common denominator $4(x+y-1)^2$, and the logarithms combine to $\frac12\log(\Theta(x)\Theta(y))$ with $\Theta(w)=(2w+1)(w+1)^2/w^2$; hence
\[
\mathcal S=\frac{t^2-10t+3s^2+1}{2(s-1)^2}-\frac{C_{\mathrm g}}2\log\frac{(4t+2s+1)(t+s+1)^2}{t^2}.
\]
At fixed $s$,
\[
\begin{aligned}
\partial_t\mathcal S&=\frac{t-5}{(s-1)^2}-\frac{C_{\mathrm g}}2\Bigl(\frac4{4t+2s+1}+\frac2{t+s+1}-\frac2t\Bigr)\\
&<-\frac{5-t}{(s-1)^2}+\frac{C_{\mathrm g}(s+1)}{t(t+s+1)},
\end{aligned}
\]
which is negative as soon as $\mathcal K(s,t)=(5-t)t(t+s+1)/((s+1)(s-1)^2)\ge\frac45>C_{\mathrm g}$. On $[1,2]^2$: if $2<s\le3$ then $s-1\le t\le9/4$ and $\mathcal K\ge\frac{11}4\cdot\frac{2s(s-1)}{(s+1)(s-1)^2}\ge\frac{11}8$; if $3\le s\le4$ then $2s-4\le t\le s^2/4$, $t+s+1\ge3(s-1)$, and $\mathcal K\ge3(5-s^2/4)(2s-4)/(s^2-1)\ge\frac45$, because $15(20-s^2)(2s-4)-16(s^2-1)=(4-s)(30s^2+76s-296)\ge0$. So at fixed $s$ the minimum is on the diagonal $x=y=w$, where
\begin{gather*}
\mathcal S(w,w)=\frac{(w^2+1)^2}{2(2w-1)^2}-C_{\mathrm g}\log\Theta(w),\\
\frac{d}{dw}\mathcal S(w,w)=2(w^2-w-1)\Bigl(\frac{w^2+1}{(2w-1)^3}-\frac{C_{\mathrm g}}{w(w+1)(2w+1)}\Bigr).
\end{gather*}
The bracket is positive on $[1,2]$, since $2w-1\le3w/2$, $w+1\ge w$, $2w+1\ge2w$, and $w^2+1\ge2w$ give $w(w+1)(2w+1)(w^2+1)/(2w-1)^3\ge32w/27>C_{\mathrm g}$. So the diagonal has its unique minimum at $w=\tau$, where $(\tau^2+1)^2=5\tau^2$, $(2\tau-1)^2=5$, and $\Theta(\tau)=\tau^5$, so $\mathcal S(\tau,\tau)=\tau^2/2-5C_{\mathrm g}\log\tau=0$.
\end{proof}

\subsection{Proofs of Lemma~\ref{lem:noflip} and Theorem~\ref{thm:sparse}}\label{app:sparse}

\emph{Lemma~\ref{lem:noflip}.} A no-flip state has $x\ge2$. Squaring $2.236068$ shows $\sqrt5<2.236068$, and $\log\tau>0.4812$ (proved by exact rational arithmetic in \ref{app:golden}); hence $C_{\mathrm g}<5.236068/9.624<\kappa-1/30000$, an exact comparison of rationals. Also $C_{\mathrm g}<3/4$. By Cauchy--Schwarz $\Gamma\ge\frac18$, since $2(x^2+y^2)\ge(x+y)^2>(x+y-1)^2$. If $\Lambda(x,y)\le\frac1{12}$ then $\varrho\ge\frac18-\frac34\cdot\frac1{12}=\frac1{16}$. Otherwise Lemma~\ref{lem:outside} gives $\varrho>(\kappa-C_{\mathrm g})\Lambda>\frac1{30000}\cdot\frac1{12}=\frac1{360000}$. The bound on $\Phi$ follows from Proposition~\ref{prop:action} and Lemma~\ref{lem:discharge}, since all no-flip residuals are among the nonnegative ones.

\emph{Theorem~\ref{thm:sparse}.} Let $p$ qualify and let $r$ attain $m(p)$; by Proposition~\ref{prop:sym} we may take $r\le p/2$. Write $p/r=[a_1;\dots,a_K]$ canonically and $E=\sum_i(a_i-1)$, so $N_0=E-1$ by Proposition~\ref{prop:action}. Lemma~\ref{lem:noflip} gives $(E-1)/1080000\le\varepsilon_0\log p+C_{\mathrm g}\log2-\frac18$, that is,
\[
E\le\frac1{10}\log p+D_1,\qquad D_1=1080000\Bigl(C_{\mathrm g}\log2-\frac18\Bigr)+1>0 .
\]
By \eqref{eq:length} and $\log\tau>1/3$, $K<3\log p$. Put $L=3\log X+1$ and $E_X=\frac1{10}\log X+D_1$. Every qualifying $p\le X$ is the numerator of a word $[a_1;\dots,a_k]$ of positive integers with $k\le L$ and $\sum_i(a_i-1)\le E_X$, and a word determines its value and hence $p$. Writing $a_i=1+b_i$ and adding a slack variable, there are $\binom{k+\lfloor E_X\rfloor}{\lfloor E_X\rfloor}$ such words of length $k$, so for $X$ large the number of qualifying $p$ is at most
\[
L\binom{\lfloor L\rfloor+\lfloor E_X\rfloor}{\lfloor E_X\rfloor}\le L\Bigl(\frac{e(L+E_X)}{E_X}\Bigr)^{E_X},
\]
because $\binom ab\le(ea/b)^b$ and $y\mapsto y\log(e(L+y)/y)$ is increasing. Here $(L+E_X)/E_X\le31+10/\log X$. Splitting the exponent $E_X=\frac1{10}\log X+D_1$, and using $(1+\frac{10}{31\log X})^{\frac1{10}\log X}\le e^{1/31}$ on the first part and $31+10/\log X\le32$ on the second, gives $(e(L+E_X)/E_X)^{E_X}\le(31e)^{\frac1{10}\log X}\,e^{1/31}\,(32e)^{D_1}$ for $X\ge e^{10}$. The count is therefore $O(X^\theta\log X)$ with $\theta=(1+\log31)/10$, and $\theta<\frac12$ because the partial sum $\sum_{j=0}^44^j/j!=103/3$ exceeds $31$, so $\log31<4$. The remaining assertions follow, using Theorem~\ref{thm:golden} for the liminf.
\qed

\subsection{Proof of Proposition~\ref{prop:minimizers}}\label{app:minimizers}

The digit-sum bound is the lower bound of Proposition~\ref{prop:order}. In the chain of $(n_0,n_1)=(p,r)$ consider block $i$ with $a=a_i$, $n=n_i$, $s=n_{i+1}<n$, $u=u_i$, and $v=u_{i-1}\le u$. The identity $n_{i-1}u_i+n_iu_{i-1}=p$, which holds for $i=1$ and is preserved by the recursions, gives $p=u(an+s)+vn\le un(a+2)$. By Theorem~\ref{thm:closed-form} (with $(q,p)$ replaced by $(p,r)$), $\Phi(p,r)\ge u^2P(n,s,a)\ge u^2\sum_{j=1}^a\binom{jn}2\ge u^2n^2a(a^2-1)/6$, so $\Phi(p,r)\ge p^2a(a^2-1)/(6(a+2)^2)$, and $a(a^2-1)/(a+2)^2=a-4+(11a+16)/(a+2)^2>a-4$. For the first block $v=0$ and for the last $s=0$, so $p\le un(a+1)$ and $\Phi(p,r)\ge p^2a(a-1)/(6(a+1))>p^2(a-2)/6$. Finally $p/r<a_1+1<6H/p^2+3$, and $p-r$, $r^*$, $p-r^*$ satisfy the same bounds by Proposition~\ref{prop:sym}. \qed

\subsection{Proof of Theorem~\ref{thm:penalty}}\label{app:penalty}

\emph{Certificate:} \texttt{certificates/mixed\_stencil}.

For a step with neighbors, call the incoming move and the outgoing move flips (F) or no-flips (N); the outgoing move is N exactly when $x\ge2$, and for an interior step the incoming move is N exactly when $y\ge2$, its predecessor being $\psi_{\mathrm n}^{-1}(x,y)=(x+1,y-1)$ or $\psi_{\mathrm f}^{-1}(x,y)$. Let $\mathcal S_j=\varrho_{j-1}+3\varrho_j+\varrho_{j+1}$.

\emph{Three stencil bounds.} (FF) $\mathcal S_j\ge0$ by Lemma~\ref{lem:stencil}. (NN) For $x,y\ge2$, $\mathcal S_{NN}(x,y)=\varrho(x+1,y-1)+3\varrho(x,y)+\varrho(x-1,y+1)>1/40$. Indeed, with $f_N(w)=(w+2)(w+1)^2/(w^2(w-1))$,
\begin{gather*}
\mathcal S_{NN}=\frac{5(x^2+y^2)+4}{4(x+y-1)^2}-\frac{C_{\mathrm g}}2\log\bigl(f_N(x)f_N(y)\bigr),\\
f_N(x)f_N(y)=\frac{(t+2s+4)(t+s+1)^2}{t^2(t-s+1)},
\end{gather*}
with $s=x+y$, $t=xy\in[2s-4,s^2/4]$. The rational part is affine in $t$, and $\partial_t^2\mathcal S_{NN}=\frac{C_{\mathrm g}}2\bigl[(t+2s+4)^{-2}+2(t+s+1)^{-2}-2t^{-2}-(t-s+1)^{-2}\bigr]<0$, so the minimum at fixed $s$ is at $x=y$ or at $x=2$. With $h(w)=(5w^2+5w-4)/(w(w+1)(w-1)(w+2))>0$, the diagonal function $\frac{5w^2+2}{2(2w-1)^2}-C_{\mathrm g}\log f_N(w)$ has derivative $-\frac{5w+4}{(2w-1)^3}+C_{\mathrm g}h(w)$, and the boundary function $\frac{5w^2+24}{4(w+1)^2}-\frac{C_{\mathrm g}}2\log(9f_N(w))$ has derivative $\frac{5w-24}{2(w+1)^3}+\frac{C_{\mathrm g}}2h(w)$. Using $C_{\mathrm g}>\frac12$ (from $\sqrt5<\frac94$, $\tau<\frac{13}8$, and $\log\frac{13}8\le2a+\frac{2a^3}{3(1-a^2)}<\frac12$ with $a=\frac5{21}$, together with $\tau^2=\tau+1>\frac52$), both derivatives are positive for $w\ge2$, because
\begin{align*}
&(5w^2+5w-4)(2w-1)^3-2(5w+4)w(w+1)(w-1)(w+2)\\
&\qquad=30+439v+941v^2+748v^3+252v^4+30v^5,\\
&(5w^2+5w-4)(w+1)^2+2(5w-24)w(w-1)(w+2)\\
&\qquad=10+69v+165v^2+97v^3+15v^4,
\end{align*}
with $w=2+v$. Hence $\mathcal S_{NN}\ge\mathcal S_{NN}(2,2)=\frac{11}9-C_{\mathrm g}\log9>\frac{11}9-\kappa\cdot\frac{11}5=\frac{11341}{450000}>\frac1{40}$, using $\log3<2(\frac12+\frac1{24})+\frac{2(1/2)^5}{5(3/4)}=\frac{11}{10}$.

(FN and NF) For $x\ge2$, $1<y<2$, $\mathcal S_{FN}(x,y)=\varrho(\psi_{\mathrm f}^{-1}(x,y))+3\varrho(x,y)+\varrho(x-1,y+1)>1/100$, and $\mathcal S_{NF}(x,y)=\mathcal S_{FN}(y,x)$ because exchanging the coordinates conjugates $\psi_{\mathrm f}$ to $\psi_{\mathrm f}^{-1}$ and $\psi_{\mathrm n}$ to $\psi_{\mathrm n}^{-1}$. This bound is computer-assisted. With $t=1/x\in(0,1/2]$ and $z=y-1\in(0,1)$, the rational part of $\mathcal S_{FN}$ is $n(t,z)/(4(1+tz)^2)$ with $n=z^2+5+2t(z^2-1)+t^2(5z^2+10z+8)$, and its ten logarithms have the six arguments $1+t$, $1/t$ (three times), $(1-t)/t$, $1/z$, $1+z$ (three times), and $2+z$, all at least $1$. Replacing each $\log(1+1/u)$ by $f(u)$ and $C_{\mathrm g}$ by $\kappa$ gives $\mathcal S_{FN}\ge\frac{n}{4(1+tz)^2}-\frac\kappa{12}(a(t)+b(z))$ with
\begin{align*}
a(t)&=4t+\frac{1+3t}{1+t}+\frac{1+24t}{2+t}+\frac8{3+2t}+\frac t{1-t}+\frac{8t}{2-t},\\
b(z)&=z+\frac{z+3}{1+z}+\frac{8z+4}{2+z}+\frac{24}{3+2z}+\frac1{3+z}+\frac8{5+2z}.
\end{align*}
Let $D_t=(1+t)(2+t)(3+2t)(1-t)(2-t)$, $D_z=(1+z)(2+z)(3+2z)(3+z)(5+2z)$, $D=D_tD_z$, $a_N=D_ta$, and $b_N=D_zb$. Clearing the positive factor $120000D(1+tz)^2$, it suffices that
\[
\begin{split}\mathcal P&=30000Dn-1200D(1+tz)^2-5441(1+tz)^2(a_ND_z+b_ND_t)\\&>0\qquad\text{on }[0,\tfrac12]\times[0,1].\end{split}
\]
$\mathcal P$ has degree eight in each variable. On each of the sixteen boxes $[\frac i8,\frac{i+1}8]\times[\frac j4,\frac{j+1}4]$ the $81$ tensor Bernstein coefficients of degree $(8,8)$ were computed in exact rational arithmetic by two independent implementations; all $1296$ are positive (the least is $89018097/32768$), which proves $\mathcal P>0$.

\emph{Assembly.} Let the path have $n\ge1$ steps, let $Q$ be the number of no-flip steps (including the last step, from $x=2$ to the terminal state), and add the virtual states $v_-=(x_1+1,1)$ and $v_+=(1,y_n+1)$ as neighbors of the first and last steps; these are $\psi_{\mathrm n}^{-1}(x_1,y_1)$ and $\psi_{\mathrm n}(x_n,y_n)$, since $y_1=2$ and $x_n=2$. By the three bounds, $\mathcal S_j\ge\frac1{100}(\mathbf 1_{\text{in N}}+\mathbf 1_{\text{out N}})$ at every step, and the total number of N incidences is $2Q$ ($Q$ outgoing, $Q-1$ internal incoming plus the virtual one). Hence $\sum_j\mathcal S_j\ge Q/50$. Also $\sum_j\mathcal S_j=5\sum_j\varrho_j-\varrho_1-\varrho_n+\varrho(v_-)+\varrho(v_+)$, $\varrho_1,\varrho_n\ge0$ by Lemma~\ref{lem:outside}, and for $w\ge2$, with $u=1/w$,
\[
\varrho(w,1)=\varrho(1,w)=\tfrac14+\tfrac{u^2}4-\tfrac{C_{\mathrm g}}2(\log2+\log(1+u))\le\tfrac14-\tfrac{C_{\mathrm g}}2\log2=:J_0,
\]
because $\log(1+u)\ge u/(1+u)$ and $u(1+u)\le\frac34<2C_{\mathrm g}$. Therefore $\sum_j\varrho_j\ge Q/250-2J_0/5$, and Proposition~\ref{prop:action} gives
\[
\frac{\Phi(p,r)}{p^2}\ge C_{\mathrm g}\log\frac p2+\frac Q{250}-\frac{2J_0}5+\frac18+\frac1{4p}=C_{\mathrm g}\log p+\frac Q{250}-D_0+\frac1{4p}.
\]
For $p=2$ ($Q=0$) the inequality is checked directly. By Proposition~\ref{prop:action}, $Q=E-1$. Finally $D_0<\frac45\kappa\cdot\frac{1733}{2500}-\frac1{40}=\frac{8648003}{31250000}<\frac{277}{1000}$, using $\log2\le2(a+\frac{a^3}3+\frac{a^5}5)+\frac{2a^7}{7(1-a^2)}<\frac{1733}{2500}$ with $a=\frac13$.

\emph{Counting.} If $m(p)\le(C_{\mathrm g}+\varepsilon)p^2\log p$ and $r\le p/2$ is a minimizer, then $E<250\varepsilon\log p+\frac{281}4$. Let $\mathsf A=1/\log\tau<\frac{25}{12}$. Every such $p\le X$ is the numerator of a word of positive integers of length $k\le\mathsf A\log X$ and excess $\sum(a_i-1)\le E_X:=\beta\log X+\frac{281}4$, $\beta=250\varepsilon$. For $0<z<1$ the number of words of length $k$ and excess at most $E_X$ is at most $z^{-E_X}(1-z)^{-k}$, so, summing the geometric series over $k$, the number of such $p$ is at most $z^{-285/4}X^{\beta\log(1/z)+\mathsf A\log(1/(1-z))}$. For $\varepsilon=\frac1{200000}$ take $z=\frac1{1000}$: since $\log1000<7$ (as $\sum_{j\le11}7^j/j!>1000$) and $\log\frac{1000}{999}<\frac1{999}$, the exponent is less than $\frac7{800}+\frac{25}{12\cdot999}<\frac1{90}$. For $\varepsilon=\frac1{750}$ take $z=\frac17$: since $\log7<\frac{39}{20}$ (as $\sum_{j\le9}(39/20)^j/j!>7$) and $\log\frac76=2\operatorname{artanh}\frac1{13}<\frac{31}{200}$, the exponent is less than $\frac{13}{20}+\frac{31}{96}<\frac{39}{40}$. For $\varepsilon=0$, $E\le70$. A word of length one is $p/r=[p]$ with $r=1$ and $E=p-1$, so it gives $p\le71$. In a word of length $k\ge2$, $a_1\ge2$ because $r<p/2$ (as $p>2$), and $a_k\ge2$ because the expansion is canonical. Writing $a_1=2+b_1$, $a_k=2+b_k$, and $a_i=1+b_i$ for $1<i<k$, the integers $b_i\ge0$ have sum $E-2\le68$, and with a slack variable there are $\binom{k+68}{68}$ such words. Hence the number of these $p\le X$ is at most $70+\sum_{2\le k\le\mathsf A\log X}\binom{k+68}{68}=O((\log X)^{69})$. The exceptional counts are $o(X/\log X)$, and $\pi(X)\gg X/\log X$. \qed

\subsection{Proof of Theorem~\ref{thm:typical-upper}}\label{app:typical}

\emph{Certificate:} \texttt{certificates/transfer\_73\_125}.

\emph{A cocycle formula.} Let $r/p=[0;a_1,\dots,a_k]$ be any expansion with positive digits (canonical, or ending in $1$), with continuants $q_i$, $v=q_{k-1}$, past tails $\omega_i=q_{i-2}/q_{i-1}=[0;a_{i-1},\dots,a_1]$, and future tails $\xi_i=[0;a_{i+1},\dots,a_k]$ (an empty tail is $0$). Then
\begin{equation}\label{eq:cocycle}
\frac{\Phi(p,r)}{p^2}=\sum_{i=1}^kg(a_i,\xi_i,\omega_i)-\frac{p+r+v-1}{4p},\qquad \sum_{i=1}^k\log(a_i+\omega_i)=\log p,
\end{equation}
\[
g(a,\xi,\omega)=\frac{a(a+1)(2a+1)/3+a(a+1)(\xi+\omega)+a(\xi^2+\omega^2)}{4(a+\xi+\omega)^2}.
\]
The second identity is $q_i/q_{i-1}=a_i+\omega_i$. For the first, run the subtractive algorithm from lengths $(p,r)$ and counts $(0,1)$, performing $a_i$ subtractions in block $i$ and then interchanging, and continue one step past $(1,1)$ to $(0,1)$. At the start of block $i$, $A/B=a_i+\xi_i$ and $U/V=\omega_i$, so its steps have $x=a_i+\xi_i-\ell$ and $y=\omega_i+\ell+1$ ($0\le\ell<a_i$) with constant $x+y-1=a_i+\xi_i+\omega_i$, and the sum of $\Gamma$ over the block is $g(a_i,\xi_i,\omega_i)$. Compared with the path of Proposition~\ref{prop:action}, this adds the initial virtual state $(x,y)=(p/r,1)$, with $\Gamma=(p^2+r^2)/(4p^2)$, and the final virtual state $(1,p/v)$, with $\Gamma=(p^2+v^2)/(4p^2)$, the terminal counts being $\{p-v,v\}$. Subtracting these from the boundary term of Proposition~\ref{prop:action} with $U_*V_*=(p-v)v$ gives $-(p+r+v-1)/(4p)$.

\emph{Monotonicity.} Up to a positive factor, $\partial_\xi g=(a-1+2\omega)\xi-2\omega^2-(a+1)\omega-(a+1)(a+2)/3\le-(a^2+5)/3-(a-1)\omega-2\omega^2<0$ for $\xi,\omega\in[0,1]$, and symmetrically for $\omega$. Hence $f(a,\xi,\omega)=4g(a,\xi,\omega)-\frac{489}{125}\log(a+\omega)$ decreases in both tails.

\emph{The transfer certificate.} Restrict to the alphabet $\{1,2,3,4\}$. For a block $b$ of four digits let $\ell(b)$ be the smaller of $[0;b,0]$ and $[0;b,1]$ (composition of the maps $w\mapsto1/(a+w)$), so that every tail beginning with $b$ is at least $\ell(b)$. On the $4^8$ words of eight digits define a matrix $T$ whose entry from $w$ to the word obtained by deleting the first digit of $w$ and appending a digit $d$ is $e^{\hat f}$, where $\hat f$ is $f$ evaluated at the fifth digit $a$ of the nine-digit window $wd$ with $\xi=\ell(\text{last four digits})$ and $\omega=\ell(\text{first four digits, reversed})$. Then $f\le\hat f$ at every digit of a word that has at least four digits on each side. A positive integer vector $h$ satisfies
\[
Th<\tfrac{999}{1000}\,h,
\]
which was verified row by row in exact integer arithmetic: each $\hat f$ is bounded above using rational lower bounds for logarithms from the series of $\operatorname{artanh}$, and each $e^{\hat f}$ is bounded above, after rounding $\hat f$ up to a multiple of $10^{-12}$, by a degree-$30$ Taylor polynomial at $\hat f/16$ with an explicit tail bound, evaluated with outward rounding and squared four times. The vector $h$ was found by power iteration on a floating-point copy of $T$ and scaled and rounded to integers; a script carrying out this procedure is part of the code, but the inequality is proved from the stored integer entries alone, whatever their origin. The largest of the certified upper bounds for the row ratios $(Th)_w/h_w$ is $0.99868\ldots$, so the inequality holds with a margin of about $0.0013$ below $1$, which is all the argument needs. A second, separately written program proves the inequality again on $163$ of the $65536$ rows. This is the computer-assisted step.

\emph{Counting expensive words.} For a word of length $k\ge8$, the digits in positions $5,\dots,k-4$ correspond to a path of length $k-8$ in the graph of $T$, and the eight remaining positions contribute at most $8\cdot\frac{15}4=30$ to $\sum_if$, since $g\le g(a,0,0)=(a+1)(2a+1)/(12a)\le\frac{15}{16}$ and $\log(a+\omega)\ge0$. Since $T^jh\le(\frac{999}{1000})^jh$ and $\mathbf 1\le h/\min_wh_w$, the words of length $k\ge8$ contribute at most $e^{30}(\frac{999}{1000})^{k-8}\sum_wh_w/\min_wh_w$ in total, and each of the finitely many shorter words at most $e^{15k/4}$; hence, by \eqref{eq:cocycle},
\[
\mathcal Z=\sum_{\mathbf a\in\{1,2,3,4\}^*}e^{4\sum_ig(a_i,\xi_i,\omega_i)}\,q(\mathbf a)^{-489/125}=\sum_{\mathbf a}e^{\sum_if(a_i,\xi_i,\omega_i)}<\infty ,
\]
where $q(\mathbf a)$ is the denominator of the word $\mathbf a$. For a word with value $r/p$, so $p=q(\mathbf a)$: if $\Phi(p,r)\ge\frac{73}{125}p^2\log p$, then \eqref{eq:cocycle} gives $\sum_ig>\frac{73}{125}\log p$, so the summand of $\mathbf a$ exceeds $p^{4\cdot73/125-489/125}=p^{-197/125}\ge X^{-197/125}$ when $p\le X$. Hence the words $\mathbf a$ over $\{1,2,3,4\}$ with $q(\mathbf a)\le X$ and $\Phi(p,r)\ge\frac{73}{125}p^2\log p$ number at most $\mathcal ZX^{197/125}$.

\emph{Arithmetic input.} Kan \cite[Thm.~1.4, formulas (1.6)--(1.8)]{Kan2021} proves that for a finite alphabet whose Cantor set has Hausdorff dimension $\delta\ge(\sqrt{40}-4)/3$, all but $O(Xe^{-c_0\sqrt{\log X}})$ admissible integers in $[X/2,X]$ have $\gg X^{2\delta-1.001}$ distinct numerators $r$ such that $r/p$ has an expansion with digits in the alphabet. Here (1.8) bounds the multiplicity $\kappa_{\mathcal A}(p)$, which \cite[p.~922 and (1.2)]{Kan2021} defines as the number of integers $1\le r\le p$ with $\gcd(r,p)=1$ and $r/p=[0;d_1,\dots,d_k]$, all $d_j$ in the alphabet; the abstract of \cite{Kan2021} states only the positive-proportion form (1.6), and Theorem~1.4 there gives all three formulas. For $\{1,2,3,4\}$, Pollicott and Vytnova \cite[Thm.~4.14]{PollicottVytnova2022} prove $\delta=0.788945557483\pm10^{-12}$. Every integer is admissible for $\{1,2\}$, hence for every alphabet containing it. Indeed, with $M_1=\left(\begin{smallmatrix}0&1\\1&1\end{smallmatrix}\right)$ and $M_2=\left(\begin{smallmatrix}0&1\\1&2\end{smallmatrix}\right)$, whose products along a word have the denominator as lower right entry and the numerator as upper right entry, the semigroup of words of even length over $\{1,2\}\subset\{1,2,3,4\}$ lies in $\mathrm{SL}_2(\Z)$ and contains $M_1^2$, $M_1M_2$, $M_2M_1$. The group they generate contains $M_1^{-1}M_2=M_1^{-2}(M_1M_2)=\left(\begin{smallmatrix}1&1\\0&1\end{smallmatrix}\right)=:\mathsf U$ and, since $M_1^{-1}\mathsf UM_1=M_1^{-2}(M_2M_1)$, also $(M_1^{-1}\mathsf UM_1)\mathsf U^{-2}=\left(\begin{smallmatrix}0&-1\\1&0\end{smallmatrix}\right)=:\mathsf S$; these generate $\mathrm{SL}_2(\Z)$. Reduction $\mathrm{SL}_2(\Z)\to\mathrm{SL}_2(\Z/N\Z)$ is surjective, and a nonempty subsemigroup of a finite group is a subgroup, so the even words reduce onto all of $\mathrm{SL}_2(\Z/N\Z)$. The word reducing to $\mathsf S\mathsf U^d=\left(\begin{smallmatrix}0&-1\\1&d\end{smallmatrix}\right)$ is a word over $\{1,2\}$ whose lower right entry, the denominator of a reduced fraction in $(0,1]$, is congruent to $d$ modulo $N$, for every $d$ and $N$.

If a non-exceptional $p\in[X/2,X]$ has $m(p)\ge\frac{73}{125}p^2\log p$, then each of its $\gg X^{2\delta-1.001}$ numerators gives at least one expensive word, and distinct numerators give distinct words. Hence the number of such $p$ is $O(X^{197/125+1.001-2\delta})=O(X^{1-1/1250})$, since $2\cdot0.7889-0.001-\frac{197}{125}=\frac1{1250}$ and $\delta>0.7889$. Summing over dyadic intervals gives the bound of Theorem~\ref{thm:typical-upper}, with a smaller $c_0$; it is $o(X/\log X)$, which gives the statement for primes, and the lower bound is Theorem~\ref{thm:penalty}. \qed

\section*{Funding}
This research did not receive any specific grant from funding agencies in the public, commercial, or not-for-profit sectors.

\section*{Declaration of competing interest}
The authors declare that they have no known competing financial interests or personal relationships that could have appeared to influence the work reported in this paper.

\section*{Data availability}
The code, the outputs of all computations, and the certificates of the computer-assisted results are archived on Zenodo at \url{https://doi.org/10.5281/zenodo.23194312} \cite{code2026}, from the tagged release of the repository \url{https://github.com/methodtrace/2026-spacings-of-multiples}.

\section*{Declaration of generative AI and AI-assisted technologies in the manuscript preparation process}
During the preparation of this work the authors used Claude Opus (Anthropic) in order to search the literature, assist in writing and checking proofs and text, and write the verification code, under the authors' direction. After using this tool, the authors reviewed and edited the content as needed and take full responsibility for the content of the published article.

\bibliographystyle{elsarticle-num}
\bibliography{refs}

\begin{thebibliography}{10}
\expandafter\ifx\csname url\endcsname\relax
  \def\url#1{\texttt{#1}}\fi
\expandafter\ifx\csname urlprefix\endcsname\relax\def\urlprefix{URL }\fi
\expandafter\ifx\csname href\endcsname\relax
  \def\href#1#2{#2} \def\path#1{#1}\fi

\bibitem{Sos1958}
V.~T. S{\'o}s, On the distribution mod 1 of the sequence $n\alpha$, Ann. Univ.
  Sci. Budapest. E{\"o}tv{\"o}s Sect. Math. 1 (1958) 127--134.

\bibitem{Suranyi1958}
J.~Sur\'anyi, {\"U}ber die {A}nordnung der {V}ielfachen einer reellen {Z}ahl
  mod 1, Ann. Univ. Sci. Budapest. E{\"o}tv{\"o}s Sect. Math. 1 (1958)
  107--111.

\bibitem{Swierczkowski1958}
S.~{\'S}wierczkowski, On successive settings of an arc on the circumference of
  a circle, Fund. Math. 46 (1958) 187--189.
\newblock \href {https://doi.org/10.4064/fm-46-2-187-189}
  {\path{doi:10.4064/fm-46-2-187-189}}.

\bibitem{Undershoot}
P.~Chaudhary, K.~Khandelwal, Undershoot of $h$-fold sumsets: {A}p\'ery spectra,
  similar ordering and the three-distance theorem, arXiv:2611.33046 [math.CO]
  (2026).

\bibitem{Greenwood1946}
M.~Greenwood, The statistical study of infectious diseases, J. Roy. Statist.
  Soc. 109~(2) (1946) 85--110.
\newblock \href {https://doi.org/10.2307/2981176} {\path{doi:10.2307/2981176}}.

\bibitem{Moran1947}
P.~A.~P. Moran, The random division of an interval, Suppl. J. Roy. Statist.
  Soc. 9~(1) (1947) 92--98.
\newblock \href {https://doi.org/10.2307/2983572} {\path{doi:10.2307/2983572}}.

\bibitem{Shkredov2026}
I.~D. Shkredov, On some results of {K}orobov and {L}archer and {Z}aremba's
  conjecture, arXiv:2603.14116v2 [math.NT] (2026).

\bibitem{Zhang2026}
X.~Zhang, Expansion in {$\mathrm{SL}_2(\mathbb{Z}/q\mathbb{Z})$} and
  {Z}aremba's conjecture, arXiv:2605.02518v2 [math.NT] (2026).

\bibitem{Zaremba1972}
S.~K. Zaremba, La m\'ethode des ``bons treillis'' pour le calcul des
  int\'egrales multiples, in: S.~K. Zaremba (Ed.), Applications of Number
  Theory to Numerical Analysis, Academic Press, New York, 1972, pp. 39--119.

\bibitem{BourgainKontorovich2014}
J.~Bourgain, A.~Kontorovich, On {Z}aremba's conjecture, Ann. of Math. (2)
  180~(1) (2014) 137--196.
\newblock \href {https://doi.org/10.4007/annals.2014.180.1.3}
  {\path{doi:10.4007/annals.2014.180.1.3}}.

\bibitem{Huang2015}
S.~Huang, An improvement to {Z}aremba's conjecture, Geom. Funct. Anal. 25~(3)
  (2015) 860--914, arXiv:1310.3772.
\newblock \href {https://doi.org/10.1007/s00039-015-0327-6}
  {\path{doi:10.1007/s00039-015-0327-6}}.

\bibitem{Kan2021}
I.~D. Kan, A strengthening of the {B}ourgain--{K}ontorovich method: three new
  theorems, Sb. Math. 212~(7) (2021) 921--964.
\newblock \href {https://doi.org/10.1070/SM9437} {\path{doi:10.1070/SM9437}}.

\bibitem{PollicottVytnova2022}
M.~Pollicott, P.~Vytnova, Hausdorff dimension estimates applied to {L}agrange
  and {M}arkov spectra, {Z}aremba theory, and limit sets of {F}uchsian groups,
  Trans. Amer. Math. Soc. Ser. B 9 (2022) 1102--1159.
\newblock \href {https://doi.org/10.1090/btran/109}
  {\path{doi:10.1090/btran/109}}.

\bibitem{BhardwajCistoDeshpande2026}
O.~P. Bhardwaj, C.~Cisto, P.~Deshpande, Structure of gaps and generalized
  {W}ilf's conjecture, arXiv:2611.20086 [math.CO] (2026).

\bibitem{RademacherGrosswald1972}
H.~Rademacher, E.~Grosswald, Dedekind Sums, no.~16 in Carus Math. Monogr.,
  Mathematical Association of America, 1972.
\newblock \href {https://doi.org/10.5948/upo9781614440161}
  {\path{doi:10.5948/upo9781614440161}}.

\bibitem{Hickerson1977}
D.~Hickerson, Continued fractions and density results for {D}edekind sums, J.
  Reine Angew. Math. 290 (1977) 113--116.
\newblock \href {https://doi.org/10.1515/crll.1977.290.113}
  {\path{doi:10.1515/crll.1977.290.113}}.

\bibitem{Wang2020}
S.~Wang, A note on the concordance invariants {Upsilon} and phi, Commun.
  Contemp. Math. 25~(2) (2023) Paper No. 2150098, numbering as in
  arXiv:2007.11511 [math.GT].
\newblock \href {https://doi.org/10.1142/S021919972150098X}
  {\path{doi:10.1142/S021919972150098X}}.

\bibitem{vanRavenstein1988}
T.~van Ravenstein, The three gap theorem ({S}teinhaus conjecture), J. Austral.
  Math. Soc. Ser. A 45~(3) (1988) 360--370.
\newblock \href {https://doi.org/10.1017/S1446788700031062}
  {\path{doi:10.1017/S1446788700031062}}.

\bibitem{AlessandriBerthe1998}
P.~Alessandri, V.~Berth{\'e}, Three distance theorems and combinatorics on
  words, Enseign. Math. (2) 44 (1998) 103--132.

\bibitem{BreneisHinrichs2020}
S.~Breneis, A.~Hinrichs, Fibonacci lattices have minimal dispersion on the
  two-dimensional torus, in: D.~Bilyk, J.~Dick, F.~Pillichshammer (Eds.),
  Discrepancy Theory, Vol.~26 of Radon Ser. Comput. Appl. Math., De Gruyter,
  Berlin, 2020, pp. 117--132.
\newblock \href {https://doi.org/10.1515/9783110652581-006}
  {\path{doi:10.1515/9783110652581-006}}.

\bibitem{vanRavenstein1989}
T.~van Ravenstein, Optimal spacing of points on a circle, Fibonacci Quart.
  27~(1) (1989) 18--24.
\newblock \href {https://doi.org/10.1080/00150517.1989.12429596}
  {\path{doi:10.1080/00150517.1989.12429596}}.

\bibitem{Dupain1979}
Y.~Dupain, Discr{\'e}pance de la suite $(\{n\alpha\})$, $\alpha=(1+\sqrt5)/2$,
  Ann. Inst. Fourier (Grenoble) 29~(1) (1979) 81--106.
\newblock \href {https://doi.org/10.5802/aif.728} {\path{doi:10.5802/aif.728}}.

\bibitem{deBruijnErdos1949}
N.~G. de~Bruijn, P.~Erd{\H{o}}s, Sequences of points on a circle, Indag. Math.
  11 (1949) 46--49, also Proc. Kon. Ned. Akad. Wetensch. \textbf{52} (1949),
  14--17.

\bibitem{Ramshaw1978}
L.~Ramshaw, On the gap structure of sequences of points on a circle, Indag.
  Math. 40~(4) (1978) 527--541.
\newblock \href {https://doi.org/10.1016/1385-7258(78)90070-7}
  {\path{doi:10.1016/1385-7258(78)90070-7}}.

\bibitem{YaoKnuth1975}
A.~C. Yao, D.~E. Knuth, Analysis of the subtractive algorithm for greatest
  common divisors, Proc. Nat. Acad. Sci. U.S.A. 72~(12) (1975) 4720--4722.
\newblock \href {https://doi.org/10.1073/pnas.72.12.4720}
  {\path{doi:10.1073/pnas.72.12.4720}}.

\bibitem{Popov1981}
V.~N. Popov, Asymptotic formula for the sum of the sums of elements of the
  continued fractions for numbers of the form $a/p$, J. Soviet Math. 17~(5)
  (1981) 2137--2147, translated from Zap. Nauchn. Sem. Leningrad. Otdel. Mat.
  Inst. Steklov. (LOMI) \textbf{91} (1979), 81--93.
\newblock \href {https://doi.org/10.1007/BF01567592}
  {\path{doi:10.1007/BF01567592}}.

\bibitem{Panov1980}
A.~A. Panov, Averages over elements of a certain class of finite continued
  fractions, Russian Math. Surveys 35~(4) (1980) 182--183.
\newblock \href {https://doi.org/10.1070/RM1980v035n04ABEH001877}
  {\path{doi:10.1070/RM1980v035n04ABEH001877}}.

\bibitem{Panov1982}
A.~A. Panov, Mean for the sum of elements over a class of finite continued
  fractions, Math. Notes 32~(5) (1982) 781--785, translated from Mat. Zametki
  \textbf{32} (1982), no.~5, 593--600.
\newblock \href {https://doi.org/10.1007/BF01358470}
  {\path{doi:10.1007/BF01358470}}.

\bibitem{Liehl1983}
B.~Liehl, {\"U}ber die {T}eilnenner endlicher {K}ettenbr{\"u}che, Arch. Math.
  (Basel) 40 (1983) 139--147.
\newblock \href {https://doi.org/10.1007/BF01192762}
  {\path{doi:10.1007/BF01192762}}.

\bibitem{Niederreiter1978}
H.~Niederreiter, Quasi-{M}onte {C}arlo methods and pseudo-random numbers, Bull.
  Amer. Math. Soc. 84~(6) (1978) 957--1041.
\newblock \href {https://doi.org/10.1090/S0002-9904-1978-14532-7}
  {\path{doi:10.1090/S0002-9904-1978-14532-7}}.

\bibitem{Selmer1977}
E.~S. Selmer, On the linear {D}iophantine problem of {F}robenius, J. Reine
  Angew. Math. 293/294 (1977) 1--17.
\newblock \href {https://doi.org/10.1515/crll.1977.293-294.1}
  {\path{doi:10.1515/crll.1977.293-294.1}}.

\bibitem{BrownShiue1993}
T.~C. Brown, P.~J.-S. Shiue, A remark related to the {F}robenius problem,
  Fibonacci Quart. 31~(1) (1993) 32--36.
\newblock \href {https://doi.org/10.1080/00150517.1993.12429318}
  {\path{doi:10.1080/00150517.1993.12429318}}.

\bibitem{Rodseth1994}
{\O}.~J. R{\o}dseth, A note on {B}rown and {S}hiue's paper on a remark related
  to the {F}robenius problem, Fibonacci Quart. 32~(5) (1994) 407--408.
\newblock \href {https://doi.org/10.1080/00150517.1994.12429186}
  {\path{doi:10.1080/00150517.1994.12429186}}.

\bibitem{Tuenter2006}
H.~J.~H. Tuenter, The {F}robenius problem, sums of powers of integers, and
  recurrences for the {B}ernoulli numbers, J. Number Theory 117~(2) (2006)
  376--386.
\newblock \href {https://doi.org/10.1016/j.jnt.2005.06.015}
  {\path{doi:10.1016/j.jnt.2005.06.015}}.

\bibitem{code2026}
P.~Chaudhary, K.~Khandelwal, The golden ratio and the spacings of multiples
  modulo $p$: computational artifact, Zenodo, source repository:
  \url{https://github.com/methodtrace/2026-spacings-of-multiples} (2026).
\newblock \href {https://doi.org/10.5281/zenodo.23194312}
  {\path{doi:10.5281/zenodo.23194312}}.

\bibitem{RosalesGarciaSanchez2009}
J.~C. Rosales, P.~A. Garc{\'\i}a-S{\'a}nchez, Numerical Semigroups, Vol.~20 of
  Developments in Mathematics, Springer, New York, 2009.
\newblock \href {https://doi.org/10.1007/978-1-4419-0160-6}
  {\path{doi:10.1007/978-1-4419-0160-6}}.

\bibitem{RouneWoods2015}
B.~H. Roune, K.~Woods, The parametric {F}robenius problem, Electron. J. Combin.
  22~(2) (2015) P2.36.
\newblock \href {https://doi.org/10.37236/5112} {\path{doi:10.37236/5112}}.

\bibitem{Pods}
P.~Chaudhary, Unfilled seconds that demand cannot recover: duration menus,
  maxseq and live {CTV} ad pods, preprint (2026).

\bibitem{AistleitnerBordaHauke2024}
C.~Aistleitner, B.~Borda, M.~Hauke, On the distribution of partial quotients of
  reduced fractions with fixed denominator, Trans. Amer. Math. Soc. 377~(2)
  (2024) 1371--1408, arXiv:2210.14095.
\newblock \href {https://doi.org/10.1090/tran/9065}
  {\path{doi:10.1090/tran/9065}}.

\bibitem{HinrichsKritzingerPillichshammer2021}
A.~Hinrichs, R.~Kritzinger, F.~Pillichshammer, Extreme and periodic {$L_2$}
  discrepancy of plane point sets, Acta Arith. 199~(2) (2021) 163--198.
\newblock \href {https://doi.org/10.4064/aa200520-22-12}
  {\path{doi:10.4064/aa200520-22-12}}.

\bibitem{BilykTemlyakovYu2012}
D.~Bilyk, V.~N. Temlyakov, R.~Yu, Fibonacci sets and symmetrization in
  discrepancy theory, J. Complexity 28~(1) (2012) 18--36.
\newblock \href {https://doi.org/10.1016/j.jco.2011.07.001}
  {\path{doi:10.1016/j.jco.2011.07.001}}.

\bibitem{Larcher1986}
G.~Larcher, On the distribution of sequences connected with good lattice
  points, Monatsh. Math. 101~(2) (1986) 135--150.
\newblock \href {https://doi.org/10.1007/BF01298926}
  {\path{doi:10.1007/BF01298926}}.

\bibitem{Rukavishnikova2006}
M.~G. Rukavishnikova, A probability estimate for the sum of incomplete partial
  quotients with fixed denominator, Chebyshevski{\u\i} Sb. 7~(4) (2006)
  113--121, in Russian.

\bibitem{Rukavishnikova2011}
M.~G. Rukavishnikova, The law of large numbers for the sum of the partial
  quotients of a rational number with fixed denominator, Math. Notes 90~(3)
  (2011) 418--430, translated from Mat. Zametki \textbf{90} (2011), no.~3,
  431--444.
\newblock \href {https://doi.org/10.1134/S0001434611090100}
  {\path{doi:10.1134/S0001434611090100}}.

\end{thebibliography}

\end{document}